\documentclass[a4paper,11pt]{article}
\usepackage{amssymb,amsthm,bbm,epic,eepic,graphics}
\usepackage[francais]{babel}
\usepackage[applemac]{inputenc}
\usepackage{a4wide}
\usepackage{amssymb,amsthm}
\usepackage{xypic}

\def\stvrule{\vrule width-0.4pt}         
\def\sthrule{\hrule \hrule height-0.4pt} 

\def\ststrut{\vrule height1.6ex 
                     depth0.6ex 
                       width0ex 
\relax}
\def\scarre#1{%
    \vcenter{\hbox{}\hrule
             \hbox{\vrule\makebox[2.3ex]{\ststrut$\scriptstyle#1$}\vrule}\sthrule}%
             \stvrule}
\def\sgenruban#1{\vcenter{\halign{&$\scarre{##}$\cr#1}}\egroup}

\def\smallruban{%
  \bgroup
  \let\ =\omit
  \let\\=\cr
  \offinterlineskip
  \sgenruban}

\def\Stvrule{\vrule width-0.4pt}         
\def\Sthrule{\hrule \hrule height-1.6pt} 

\def\Ststrut{\vrule height0.9ex 
                     depth0.3ex 
                       width0ex 
\relax}
\def\Scarre#1{%
    \vcenter{\hbox{}\hrule
             \hbox{\vrule\makebox[1.0ex]{\Ststrut$\scriptstyle#1$}\vrule}\Sthrule}%
             \Stvrule}

\def\Sgenruban#1{\vcenter{\halign{&$\Scarre{##}$\cr#1}}\egroup}

\def\Smallruban{%
  \bgroup
  \let\ =\omit
  \let\\=\cr
  \offinterlineskip
  \Sgenruban}

\def\rtvrule{\vrule width-0.8pt}         
\def\rthrule{\hrule \hrule height-0.8pt} 

\def\rtstrut{\vrule height3.2ex 
                     depth1.2ex 
                       width0ex 
\relax}
\def\rcarre#1{%
    \vcenter{\hbox{}\hrule
             \hbox{\vrule\makebox[4.6ex]{\rtstrut$\scriptstyle#1$}\vrule}\rthrule}%
             \rtvrule}
\def\rgenruban#1{\vcenter{\halign{&$\rcarre{##}$\cr#1}}\egroup}

\def\ruban{%
  \bgroup
  \let\ =\omit
  \let\\=\cr
  \offinterlineskip
  \rgenruban}

\def\Lie{\mathrm{Lie}}
\def\rg{\mathrm{rg}}

\def\HRn{\mathrm{(HR_{n})}}
\def\HRm{\mathrm{(HR_{n-1})}}

\newcommand{\la}{\lambda}
\def\ne{\nearrow}
\def\ad{\mathrm{ad}}
\def\Ad{\mathrm{Ad}}
\def\aa{\mathbf{a}}

\def\cc{\mathbf{c}}
\def\aal{\mathbf{u}}
\def\ccl{\mathbf{v}}

\def\Id{\mathrm{Id}}
\def\Hom{\mathrm{Hom}}

\def\End{\mathrm{End}}
\def\Ker{\mathrm{Ker}\,}

\def\U{\mathsf{U}}

\def\SN{\ensuremath{\mathfrak{S}_n}}

\def\C{\ensuremath{\mathbbm{C}}}
\def\Q{\mathbbm{Q}}
\def\Z{\mathbbm{Z}}
\def\N{\mathbbm{N}}
\def\R{\mathbbm{R}}
\def\bar{\overline}
\def\det{\mathrm{det}}
\def\tr{\mathrm{tr}}
\def\rg{\mathrm{rg}}
\def\ALTN{\mathfrak{A}_n}
\def\g{\mathfrak{g}}
\def\PP{\mathfrak{P}}
\def\h{\mathfrak{h}}
\def\gl{\mathfrak{gl}}
\def\sl{\mathfrak{sl}}
\def\so{\mathfrak{so}}
\def\sp{\mathfrak{sp}}
\def\osp{\mathfrak{osp}}

\def\k{\mathbbm{k}}

\def\AN{\ensuremath{\mathfrak{B}_n}}
\def\FN{\mathbf{F}_n}

\def\om{\omega}

\def\eps{\epsilon}
\def\g{\mathfrak{g}}

\def\onto{\twoheadrightarrow}
\def\into{\hookrightarrow}

\def\dd{\mathfrak{d}}
\def\D{\mathtt{D}}
\def\a{\mathfrak{a}}
\def\cnp#1#2{\left(\begin{array}{c} {#1} \\ {#2} \end{array} \right)}
\newtheorem{prop}{Proposition}
\newtheorem{defi}{Définition}
\newtheorem{lemme}{Lemme}

\newtheorem*{theoA}{Théorème A}
\newtheorem*{theoB}{Théorème B}
\newtheorem*{theoC}{Théorème C}

\def\og{``}
\def\fg{''}

\def\s{\mathfrak{s}}

\def\B{\mathbf{B}}
\def\P{\mathbf{P}}
\def\PPhi{\widetilde{\Phi}}
\def\TTN{\widehat{\mathcal{T}_n}}
\def\Lie{\mathrm{Lie}}

\title{L'algèbre de Lie des transpositions}

\author{Ivan Marin}

\begin{document}

\maketitle

\noindent {\bf Résumé.} Pour tout $n \geq 3$ on obtient la décomposition en
facteurs simples de la sous-algèbre de Lie de l'algèbre de groupe
du groupe symétrique sur $n$ lettres engendrée par les transpositions. Cela
nous permet de déterminer l'enveloppe algébrique du groupe de tresses
$\B_n$ et de certains de ses sous-groupes dans les représentations de
l'algèbre de Iwahori-Hecke de type A.

\bigskip

\noindent {\bf Abstract.} For any $n \geq 3$ we obtain the decomposition
in simple factors of the Lie subalgebra of the group algebra of the
symmetric group on $n$ letters generated by the transpositions.
This enables us to determine the algebraic hull of the braid group
$\B_n$ and of several of its subgroups inside the representations of
the Iwahori-Hecke algebra of type A.

\section{Motivations}

On note $\B_n$ le groupe de tresses à $n$ brins. Des
actions de $\B_n$ se rencontrent dans les domaines
les plus variés de la géométrie et de l'algèbre. En particulier,
quatre types de représentations (dépendant d'un ou deux
paramètres génériques) sont particulièrement fréquents :
\begin{enumerate}
\item La représentation de Burau
\item Les représentations de l'algèbre de Temperley-Lieb
\item Les représentations de l'algèbre d'Iwahori-Hecke de type A
\item Les représentations de l'algèbre de Birman-Wenzl-Murakami
\end{enumerate}
Chacun de ces types de représentations est un cas particulier
des suivantes. Les représentations des algèbres semi-simples
correspondantes,
qui induisent des représentations du groupe de tresses,
sont bien connues en tant que représentations d'algèbres.
En revanche, la structure d'algèbre ne permet pas de comprendre
les décompositions en irréductibles des produits tensoriels
de ces représentations, en tant que représentations du groupe de tresses.
De même, elle ne permet pas de déterminer l'enveloppe
algébrique du groupe de tresses dans chacune de ces représentations.

Pour décomposer ces produits tensoriels, on a montré dans \cite{QUOTDEF}
qu'il suffit de savoir décomposer en facteurs simples
des algèbres de Lie \og infinitésimales \fg\ naturellement
associées à ces structures. Cette tâche a été accomplie
pour la représentation de Burau et l'algèbre de Temperley-Lieb
pour tout $n$, ainsi que pour
l'algèbre
d'Iwahori-Hecke pour $n \leq 6$. On y a de plus montré que
l'algèbre de Lie correspondant
à l'algèbre d'Iwahori-Hecke admet une description particulièrement
simple, précisément qu'il s'agit de la sous-algèbre de Lie de l'algèbre du groupe
du groupe symétrique, considérée comme algèbre de Lie pour le
crochet $[a,b] = ab-ba$, qui est engendrée par les transpositions.

C'est cette algèbre de Lie que nous décomposons ici (théorème A) pour tout
entier $n$.
Nous en déduisons (théorème B) l'enveloppe algébrique de l'image du groupe
de tresses
et de plusieurs de ses sous-groupes dans chacune des
représentations irréductibles de l'algèbre d'Iwahori-Hecke de type A,
puis dans l'algèbre d'Iwahori-Hecke elle-même (theorème C).

La décomposition particulièrement simple de cette
algèbre de Lie permet d'envisager le même travail
pour l'algèbre de Birman-Wenzl-Murakami, mais également
pour les algèbres de Hecke infinitésimales associées aux autres
groupes de réflexions, qui ont été introduites dans \cite{HECKINF}.

Le théorème A, établissant la décomposition de l'algèbre de Lie
des transpositions, est énoncé en section 2. La section 3 est composée de
résultats préliminaires sur les diagrammes de Young. Le but des sections 4, 5 et 6
est de démontrer ce théorème. En section 4 nous
plongeons cette algèbre de Lie dans une algèbre de Lie
réductive explicitement décomposée. Les sections 5 et 6 sont
consacrées à la démonstration de la surjectivité de
ce plongement, ce qui conclut la démonstration du
théorème. Les sections suivantes sont consacrées à
l'application de ce théorème aux calculs d'enveloppes algébriques
(théorèmes B et C).

\section{Résultat principal}

\subsection{Notations}
On rappelle qu'une partition $\la$ de l'entier $n$ est une suite
$(\la_r)_{r \geq 1}$ d'entiers naturels presque tous
nuls, tels que $\la_r \geq \la_{r+1}$ pour tout $r \geq 1$, et dont
la somme vaut $n$. On appellera plus généralement partition
une suite
$\la = (\la_r)_{r \geq 1}$ d'entiers naturels presque tous
nuls, tels que $\la_r \geq \la_{r+1}$ pour tout $r \geq 1$. C'est
une partition d'un certain entier $n$, que l'on appelle la taille
de la partition et que l'on note $|\la|$. On utilisera la notation classique
$\la \vdash n$ pour indiquer que $\la$ est une partition de $n$,
et la notation $[\la_1,\dots,\la_r]$
pour désigner la partition $\la_1,\dots,\la_r,0,\dots$. Une extension commode de cette notation
que nous utiliserons également consiste à noter $a^b$ une juxtaposition de $b$ fois le nombre $a$. Ainsi,
$[4,2^3]= [4,2,2,2]$.

Si $\la \vdash n$ est une
partition de l'entier naturel $n$, on note $\la'$ la
partition symétrique de $\la$, définie par $\la'_r = \max \{ i \geq 1 | \la_i \geq r \}$, en convenant
que $\max \emptyset = 0$. Si $\la = \la'$, on dit que $\la$ est une partition \emph{symétrique}.

Soit $\k$ un corps de caractéristique 0.
La théorie classique des représentations du groupe
symétrique $\SN$ (cf. par exemple \cite{FH}) associe à toute partition $\la$
de $n$ un $\k \SN$-module simple (de dimension finie sur $\k$), bien déterminé à isomorphisme près. Par
exemple, à la partition $[n]$ de $n$ correspond le $\k \SN$-module
trivial de dimension 1. Plus généralement, pour tout $0 \leq  r \leq n$, les partitions
de la forme $[n-r,1^r]$, appelées des \emph{équerres}, correspondent à des représentations particulières
de $\SN$. Nous appellerons les partitions qui ne sont pas des équerres des
partitions \emph{propres}.

Par souci de lisibilité, on identifiera $\la$ à ce $\k \SN$-module simple.
On notera ainsi $\dim (\la)$
la dimension sur $\k$ de ce module, et $\sl(\la)$, $\gl(\la)$ les algèbres de Lie
spéciale linéaire et linéaire sur l'espace vectoriel sous-jacent.
On utilisera la notation $\rho_{\la}$ pour désigner la
représentation associée à $\la$, c'est-à-dire le morphisme
$\SN \to GL(\la) \subset \End_{\k}(\la)$.

Pour un entier $n$, on note $E_n$ l'ensemble des
partitions propres $\la$ de $n$ telles que $\la \neq \la'$, $F_n$
l'ensemble des partitions propres symétriques. On note $\sim$ la relation
d'équivalence sur $E_n$ qui identifie $\la$ à $\la'$. Munissant l'ensemble
des partitions de $n$ de l'ordre lexicographique, on peut alors
identifier $E_n/\sim$ à l'ensemble $\{ \la \vdash n \ | \ \dim(\la)>1 \mbox{ et } \la < \la' \}$.

\subsection{L'algèbre de Lie des transpositions}

\begin{defi}
Pour tout $n \geq 3$, on note $\g_n$ la sous-algèbre de Lie de $\k \SN$
engendrée par les transpositions, $\g'_n = [\g_n,\g_n]$ son
algèbre dérivée et
$\g_{\la}$ l'image de $\g'_n$ dans les endomorphismes de $\la$.
\end{defi}

L'algèbre de Lie $\g_n$ admet des ensembles de générateurs
plus petits. Le lemme suivant, qui découle aisément de l'identité $\ad(s) \circ \ad(s)
= 2 - 2\Ad(s)$ dans $\g_n$ pour $s$ une transposition,
où $\ad(x)(y) = [x,y]$ et $\Ad(x)(y) = xyx^{-1}$, nous sera notamment
utile dans la preuve du lemme \ref{glaosp}.
\begin{lemme} \label{gentransp} $\g_n$ est engendrée par la famille des transpositions
consécutives $(i\ i+1)$. Elle est également engendrée par la
famille des transpositions $(i\ n)$ pour $1 \leq i \leq n-1$.
\end{lemme}

La proposition suivante, établie dans \cite{THESE, QUOTDEF}, montre que l'étude de $\g_n$ et de ses
représentations se ramène
à celle de $\g'_n$ :
 
\begin{prop}
L'algèbre de Lie $\g_n$ est réductive, et son centre est
de dimension 1, engendré par la somme $T_n$ de toutes les transpositions.
En conséquence $\g_n \simeq \k \times \g'_n$, et l'image de $\g_n$ dans $\gl(\la)$ est
$\g_{\la} \subset \sl(\la)$ si $T_n$ agit par $0$, et $\k \times \g_{\la}$
sinon.
\end{prop}

Si l'on note $p : \k \SN \to Z(\k \SN)$
la projection de $\k \SN$ sur son centre, définie par
$$
p(x) = \frac{1}{n!} \sum_{s \in \SN} s x s^{-1},
$$
alors l'algèbre de Lie $\g_n$ est $p$-stable parce qu'engendrée par
une classe de conjugaison de $\SN$, et on a $\k T_n = Z(\g_n) = p(\g_n)$. Ainsi,
$\g'_n$ s'identifie au noyau de la restriction de $p$ à $\g_n$.

Nous verrons en section 4.4 comment calculer l'action de $T_n$ sur chacune
des re\-pré\-sen\-ta\-tions de $\SN$. Indépendamment de ce calcul, on peut
montrer que les restrictions à $\g'_n$ de deux représentations
non isomorphes de $\k \SN$ de sont pas isomorphes en général : 
\begin{prop} \label{equivisom}
Soient $n \geq 2$ et $R_1,R_2$ deux représentations absolument irréductibles de $\k \SN$ de dimension au moins 2. On a équivalence
\begin{itemize}
\item[(i)] $R_1$ et $R_2$ sont isomorphes en tant que représentations de $\SN$.
\item[(ii)] $R_1$ et $R_2$ sont isomorphes en tant que représentations de $\g_n$.
\item[(iii)] $R_1$ et $R_2$ sont isomorphes en tant que représentations de $\g'_n$.
\end{itemize}
\end{prop}
\begin{proof} Les implications $(i) \Rightarrow (ii) \Rightarrow (iii)$
sont immédiates puisque $\g'_n \subset \g_n \subset \k \SN$. Notons
$R_i : \k \SN \to \End(V_i)$. On a par hypothèse $\dim(V_1) > 1$. Supposons
$(iii)$, c'est-à-dire qu'il existe $P \in \Hom(V_1,V_2)$ bijectif
tel que, pour tout $x \in \g'_n$, on ait
$P R_1(x) = R_2(x) P$.

A toute transposition $s$ on associe $s' = s - 2 T_n/n(n-1)$. Comme
les $n(n-1)/2$ transpositions de $\SN$ forment une classe
de conjugaison, on a $p(s) = 2 T_n/n(n-1)$ et $p(T) = T$
donc $p(s')=0$ soit $s' \in \g'_n$. D'autre part,
$T_n$ est central dans $\k \SN$ donc $R_1(T)$ et $R_2(T)$ sont
scalaires par absolue irréductibilité de ces deux représentations. De
$R_2(x) = P R_1(x) P^{-1}$ pour tout $x \in \g'_n$ on déduit
en particulier que, pour toute transposition $s$, on a
$$
R_2(s) - \frac{2}{n(n-1)} R_2(T) = PR_1(s)P^{-1} - \frac{2}{n(n-1)} R_1(T)
$$
c'est-à-dire $R_2(s) = PR_1(s)P^{-1}+ \om$ avec
$\om = 2(R_2(T_n) - R_1(T_n))/n(n-1) \in \k$ indépendant de $s$.
En élevant au carré on en déduit $1 = R_2(s^2) = 1 + 2 P R_1(s) P^{-1}
+ \om^2$. Si $\om \neq 0$, alors $R_1(s)$ serait scalaire pour toute
transposition $s$, donc pour toute permutation ; ainsi
$R_1(\SN) \subset \k^{\times}$, ce qui contredirait l'irréductibilité
de l'action de $\SN$ sur $V_1$ puisque $\dim V_1 > 1$. Ainsi $\om = 0$,
c'est-à-dire $R_2(s) = P R_1(s) P^{-1}$ pour toute transposition $s$,
donc pour toute permutation et $R_1$ est bien isomorphe à
$R_2$ comme représentation de $\SN$.
\end{proof}

On notera encore $\rho_{\la}$ la représentation $\g_n \to \gl(\la)$
associée à $\la$. Par définition, sa restriction à $\g'_n$ se factorise par $\g_{\la}$.

\subsection{Enoncé du théorème A}

Nous définirons en section 4.2, pour toute partition symétrique $\la$,
une algèbre de Lie $\osp(\la)$, et un morphisme injectif
$$
\phi_n : \g_n' \into \sl_{n-1}(\k) \times \left( \prod_{\la \in E_{n}/\sim} \sl(\la) \right) \times \left( \prod_{
\la \in F_{n}} \osp(\la) \right)
$$

Nous pouvons maintenant énoncer le theorème

\begin{theoA}  
Pour tout $n \geq 3$, $\phi_n$ est surjectif. En particulier,
$$
\g'_n \simeq \sl_{n-1}(\k) \times \left( \prod_{\la \in E_{n}/\sim} \sl(\la) \right) \times \left( \prod_{
\la \in F_{n}} \osp(\la) \right)
$$
et les représentations $\rho_{\la}$ de $\g'_n$ sont deux à deux
non isomorphes. 
\end{theoA}

\section{Préliminaires sur les diagrammes de Young}

On représente habituellement les partitions par des diagrammes de Young.
On utilisera la convention telle
que $[3,2]$ est associée au diagramme à deux colonnes et trois lignes suivante :
$$
\smallruban{ & \ \\ & \\ & \\}
$$
Si $\la$ est une partition, on notera encore $\la$ le diagramme de Young associé.

Si $\la, \mu$ sont deux partitions, on utilise les notations
suvantes :
\begin{itemize}
\item $\mu \subset \la$ si $\forall i \ \mu_i \leq \la_i$.
\item $\mu \nearrow \la$ si $\mu \subset \la$
et $|\la| = |\mu| +1$.
\item On note $\nu = \mu \cup \la$ la partition
définie par $\nu_i = \max(\la_i,\mu_i)$.
\item On note $\nu = \mu \cap \la$
la partition
définie par $\nu_i = \min(\la_i,\mu_i)$.
\end{itemize}

Comme exemple de ces deux dernières notations :
$$
\la = \smallruban{ & \ & \ & \ \\ & \ & \ & \ \\  & & \ & \ \\  & &  \ & \ \\
 & &  &  \\}\ \ \ \ 
\mu = \smallruban{ \ & \  & \ & \ \\ & & \ & \ \\  & & \ & \ \\  & &  & \ \\
 & &  & \ \\}\ \ \ \  
\la \cup \mu = \smallruban{ & \ & \ & \ \\ & & \ & \ \\  & & \ & \ \\  & &  & \ \\
 & &  &  \\}\ \ \ \  
\la \cap \mu = \smallruban{ \  & \ & \ & \ \\ & \ & \ & \ \\  & & \ & \ \\  & & \ & \ \\
 & &  & \ \\} 
$$
L'ensemble des \emph{décrochements} de $\la$ est $D(\la) = \{ i| \la_{i+1}
< \la_i \}$. On note $\delta(\la) = \# D(\la)$ le nombre de décrochements
de $\la$. A tout décrochement $r$ de $\la \vdash n$ on peut associer
$\la^{(r)} \ne \la$ définie par $\la^{(r)}_r = \la_r - 1$
et $\la^{(r)}_i = \la_i$ si $i\neq r$. Par exemple,
$$
\la = \smallruban{ & \ & \ & \ \\ & \ & \ & \ \\  & & \ & \ \\  & &  \ & \ \\
 & &  &  \\}\ \ \ 
D(\la) = \{1,2,4 \}\ \ \ 
\la^{(1)} = \smallruban{ \ & \ & \ & \ \\ & \ & \ & \ \\  & & \ & \ \\  & &  \ & \ \\
 & &  &  \\} \ 
\la^{(2)} = \smallruban{ & \ & \ & \ \\ & \ & \ & \ \\  & \ & \ & \ \\  & &  \ & \ \\
 & &  &  \\} \ 
\la^{(4)} = \smallruban{ & \ & \ & \ \\ & \ & \ & \ \\  & & \ & \ \\  & &  \ & \ \\
 & &  &  \ \\}
$$

On appelle longueur de la diagonale de $\la$ et on note $b(\la)$
le nombre $b(\la) = \max \{ i | \la_i \geq i\}$. Dans l'exemple précédent,
$b(\la) = 2$. Il est clair
que, si $\mu \subset \la$, on a $b(\mu) \leq b(\la)$, et que si $\mu \ne \la$, on
a $b(\mu) \in \{ b(\la),b(\la)-1 \}$ avec $b(\mu) = b(\la) -1$ ssi
$\mu = \la^{(b(\la))}$. On remarque que $b(\la) = 1$ signifie que
$\la = [n-r,1^r]$ pour
un certain $0 \leq r \leq n$, c'est-à-dire que $\la$ est une
\emph{équerre}.
On rappelle que dans le cas contraire on dit que $\la$ est propre,
que $E_n$ désigne l'ensemble des partitions propres non symétriques
($\la \neq \la'$) et $F_n$ celui des partitions propres symétriques.

\bigskip

Le lemme suivant est visuellement évident :

\begin{lemme} L'application $r \mapsto \la_r$ induit une bijection
$D(\la) \to D(\la')$. En particulier, $\delta(\la) = \delta(\la')$. Si
$r \in D(\la)$ et $c = \la_r$, on a $(\la')^{(c)} = (\la^{(r)})'$. Enfin
$b(\la) = b(\la')$.
\end{lemme}
\begin{proof}
Soit $r \in D(\la)$, $c = \la_r$. On a $\la'_c = \# \{ i | \la_i \geq c \}
= \# \{ i | \la_i \geq \la_r \} = r$
car $\la_r > \la_{r+1}$, et $\la'_{c+1} =
\# \{ i | \la_i \geq c+1 \}= \# \{ i | \la_i > \la_r \} \leq r-1$,
donc $\la'_{c+1} > \la'_c$ et $c \in D(\la')$. D'autre part $\la'_c = r$
s'écrit $r = \la'_{\la_r}$ donc la composée $D(\la) \to D(\la') \to
D(\la'') = D(\la)$ est l'identité et $D(\la) \to D(\la')$ est bijective.

Soit $\mu = \la^{(r)}$. Si $s \neq c$, l'inégalité $\mu_i \geq s$ signifie $\la_i \geq s$
pour tout $i \neq r$. Si $i = r$ on a $\mu_r = c-1$ et $\la_r = c$. Pour
$s > c$ on a $\mu_r < s$ et $\la_r < s$, et pour $s<c-1$
on a $\mu_r \geq c-1$ et $\la_r \geq c-1$ donc
$\mu'_i = \# \{ i | \mu_i \geq s \} = \# \{i | \la_i \geq s \} = \la'_i$.
Si $s=c$, on a $\mu'_c = \# \{ i | \mu_i \geq c \} = r-1 = \la'_c - 1$,
donc $(\la^{(r)})' = \mu' = (\la')^{(c)}$.

Enfin, $b(\la)$ est tel que $\la_{b(\la)} \geq b(\la)$. Alors
$\la'_{b(\la)} = \# \{ i | \la_i \geq b(\la) \} \geq b(\la)$.
On en déduit $b(\la') \geq b(\la)$ d'où $b(\la) = b(\la')$
puisque $\la'' = \la$.
\end{proof}

On note $P(\la) = \{ \la^{(r)} | r \in D(\la) \} = \{ \nu | \nu \ne \la \}$.
On a $\# P(\la) = \delta(\la)$. On déduit du lemme précédent
\begin{lemme} \label{lemmesym}
Si $\la = \la'$, pour tout $\mu \in P(\la)$ on a
$\mu' \in P(\la)$.
\end{lemme}
\begin{proof}
Si $\mu \in P(\la)$, c'est-à-dire s'il existe $r \in D(\la)$
tel que
$\mu = \la^{(r)}$, alors $\mu' = (\la^{(r)})' = (\la')^{\la_r}
\in P(\la') = P(\la)$.
\end{proof}

\begin{figure}
\begin{center}
\includegraphics{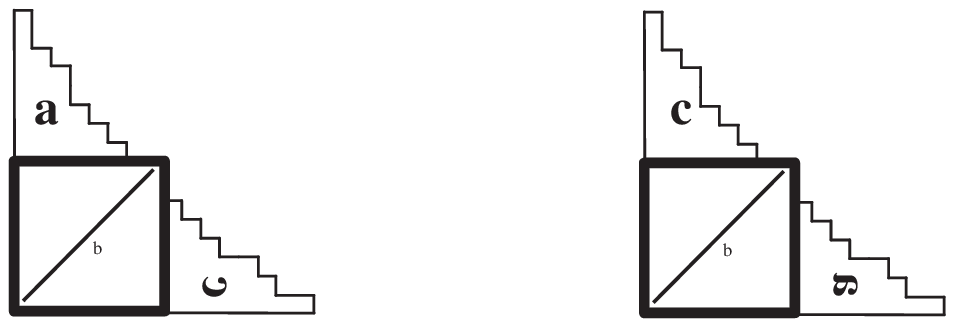}
\end{center}
\caption{Les diagrammes $\la = [\aa,\cc,b]$ et $\la' = [\cc,\aa,b]$}
\label{figacb}
\end{figure}

Dans le lemme qui suit sont regroupés un certain nombre
d'autres résultats combinatoires utiles.
Pour sa démonstration on utilisera la notation suivante.
Si $\la \vdash n$, on notera
$\la = [\aa,\cc,b]$ avec $b = b(\la)$, $\aa$ et $\cc$ deux diagrammes
de Young avec $|\la| = |\aa | + |\cc| + b^2$,
$\aa$ et $\cc$ étant déterminés par
$\la_i = b(\la) + \aa_i$ et $\la'_i = b(\la) + \cc_i$
pour $i \leq b(\la)$. On a $\la' = [\cc,\aa,b]$ (cf. figure \ref{figacb}).

\begin{lemme} \label{multilemme}
\begin{enumerate}
\item Si $b(\la) = 2$, alors $\# \{\mu \in P(\la) | b(\mu) < 2 \} \leq 1$
\item Si $\la \in F_n$, alors $P(\la) \cap F_n = \{ \la^{(b(\la))} \}$
si $b(\la) \in D(\la)$, et est vide sinon.
\item Si $\la \neq \la'$ et $\mu \neq \mu'$ avec $\mu \in P(\la)$, alors
$\mu' \not\in P(\la)$
\item $\# P(\la) \cap F_n \in \{0,1 \}$
\end{enumerate}
\end{lemme}
\begin{proof}
\begin{enumerate}
\item Si $b(\la) = 2$ et $\mu \ne \la$ tel que $b(\mu) = 1$, alors
$\mu = \la^{(b(\la))}$.
\item Si $\mu \ne \la$ et $b(\mu) < b(\la)$, on a $\mu = \la^{(b(\la))}$.
Supposons au contraire
$b(\mu) = b(\la) = b$, et écrivons $\la = [\aa,\aa,b]$ et $\mu = [\aal,\aal,
b]$. On a alors $|\aal| = |\aa|-1$ d'où $|\mu| = |\la| - 2$
ce qui est exclu.
\item  Soient $\la \neq \la'$, $\mu \neq \mu'$, $\mu \in P(\la)$ et supposons
par l'absurde $\mu' \in P(\la)$. Si $b(\mu) < b(\la)$,
$\mu = \la^{(b(\la))}$ et $\la_{b(\la)} = b(\la) = \la'_{b(\la)}$.
Posons $\mu = [\aa,\cc,b(\mu)]$. On peut supposer
$\la = [\aa^+,\cc,b(\la)]$ avec $\aa \ne \aa^+$ si $b(\mu) = b(\la)$,
et $\aa = \aa^+$ sinon. Comme $\mu' \subset \la$ on a $\cc \subset \aa^+$,
et $\aa \subset \cc$. Si $\aa = \aa^+$ on en déduit $\aa=\cc$ donc $\mu = \mu'$,
ce qui est exclu. Sinon, on a $|\aa^+| = |\aa| + 1$ donc $\aa =\cc$ ce qui
est exclu pour la même
raison, ou bien $\aa^+ = \cc$ qui est encore exclu car $\la \neq \la'$.
\item Soit $\la = [\aa,\cc,b(\la)]$. Soient $\mu,\nu \in P(\la) \cap F_n$
avec $\mu \neq \nu$. Si $b(\mu) < b(\la)$, on a $\mu = \mu' \Rightarrow \la
= \la'$. Comme $\mu \neq \nu$, on a $b(\nu) = b(\la) = b$, donc
$\la = [\aa,\aa,b]$, $\nu = [\aal,\aal,b]$, $\mu = [\aa,\aa,b-1]$
avec $\aal \ne \aa$. Mais alors $|\nu| = |\la| -2$ qui est exclu.
On peut donc supposer $b(\mu) = b(\la) = b(\nu)$. Alors
$\la = [\aa,\cc,b]$, $\mu = [\aal,\aal,b]$,
$\nu = [\ccl,\ccl,b]$. Comme $\mu \ne \la$, on peut supposer
$\aal \ne \aa$, $\aal = \cc$. Si $\ccl \ne \aa$ et $\ccl = \cc$ on aurait
$\aal = \ccl$ et $\mu = \nu$ ce qui est exclu. Comme $\nu \ne
\la$ on a donc $\ccl \ne \cc$ et $\ccl = \aa$. Mais alors
$|\cc| = |\ccl|+1 = |\aa|+1$ et $|\aa| = |\aal|+1 = |\cc|+1$, une
contradiction.
\end{enumerate}
\end{proof}

\section{Préliminaires sur les représentations de $\SN$}

Soit $\k$ un corps de caractéristique 0.
A toute partition $\la$ de l'entier $n \geq 1$ est associée
classiquement une représentation irréductible sur $\k$
du groupe symétrique $\SN$. Une description matricielle
de ces représentations est donnée par la combinatoire
des \emph{tableaux de Young standards}, que nous rappelons ici.
Pour les résultats classiques utilisés sans référence dans cette section,
on pourra consulter \cite{FH}, ch. 4.

\subsection{Modèles matriciels}

Un tableau de Young de forme $\la$ est un remplissage des cases du diagramme
de Young associé à $\la$ par les $n$ nombres de $1$ à $n$ (qui apparaissent
donc chacun une seule fois). Il est dit
\emph{standard} si la répartition de ces nombres est croissante suivant
les lignes et les colonnes. On remarque que si $T$ est un tableau
(standard) de forme $\la$, il définit de façon naturelle
un tableau de Young (standard), noté $T'$, de forme $\la'$.
Suivant notre convention, les tableaux standards associés à la
partition $[2,1]$ sont les suivants
$$
\smallruban{ 2 & \ \\ 1 & 3 \\} \ \ \ \ \ \  \smallruban{ 3 & \ \\ 1 & 2 \\}
$$
et le 2 du premier diagramme se trouve en première colonne et deuxième
ligne.

 Soit $T$ un tableau de Young de forme $\la$. Pour tout
$1 \leq r \leq n$, on note $c_r(T)$ (resp. $l_r(T)$) la colonne
(resp. la ligne)
de $T$ où se trouve $r$. Pour tous $1 \leq i<j \leq n$ on
définit la \emph{distance axiale}
$$
d_T(i,j) = c_i(T)-c_j(T) + l_j() - l_i(T)
$$
et on remarque que $d_T(j,i) = -d_T(i,j)$. De plus, lorsque $T$
est standard, si $i<j$ et $c_i(T) > c_j(T)$, on a nécessairement
$l_i(T)<l_j(T)$ donc $d_T(i,j) > 0$.

La représentation associée à $\la$
est définie sur l'espace vectoriel abstrait de base l'ensemble
des tableaux standards de forme $\la$ par l'action des transpositions
consécutives $s_r = (r\ r\! +\! 1)$. Pour ce faire, on
a besoin de la notation suivante.

Soit $T$ un tableau de Young standard de forme $\la$ et
$T_r$ le tableau déduit de $T$
en intervertissant $r$ et $r+1$. Supposons que $r$ et $r+1$ ne sont ni
dans la même ligne ni dans la même colonne. Alors $T_r$ est également
un tableau standard, et $(T_r)' = (T')_r$.
Si $c_r(T)<c_{r+1}(T)$, c'est-à-dire que $d_T(r+1,r) > 0$,
on note $T < T_r$, et $T > T_r$ dans le cas contraire.

Alors, l'action de $\SN$ dans
cette base des tableaux standards est définie
par $s_r . T = T$ si $r$ et $r+1$ sont
sur la même colonne de $T$, par $s_r . T = -T$ si $r$ et $r+1$ sont
sur la même ligne de $T$. Dans les autres cas, quitte à intervertir
$T$ et $T_r$ on peut supposer $T<T_r$ ; alors le plan engendré
par $T$ et $T_r$ est stable par l'action de $s_r$, qui est donnée dans
la base $(T,T_r)$ par la matrice
$$
\frac{1}{d} \left( \begin{array}{cc} -1 & d+1 \\ d-1 & 1 \end{array} \right)
$$ 
où $d = d_T(r+1,r)>0$.

A partir de maintenant, on fera l'abus de notation commode
qui consiste à identifier la partition ou le diagramme
de Young $\la$ avec la représentation associée, quand il n'y a
pas de risque de confusion. L'essentiel de la
théorie des représentations en caractéristique 0 du groupe
symétrique est résumé dans la proposition suivante :
\begin{prop} Soit $n \geq 1$. Pour toute partition $\la$
de $n$, la représentation de $\SN$ qui lui est associée est
absolument irréductible, et toute représentation irréductible
de $\SN$ est isomorphe à une et à une seule
d'entre elles.
\end{prop}
En particulier, toutes les représentations irréductibles
étant réalisables sur $\Q$ donc sur $\R$, elles sont toutes
autoduales.

Les propriétés importantes qui suivent nous seront utiles ici :

\paragraph{Modèle orthogonal.}
A tout tableau de Young
standard $T$ on associe l'élément suivant
$$
\zeta(T) = \prod_{\stackrel{i<j}{c_i(T)>c_j(T)}}
\frac{ d_T(i,j)-1 }{d_T(i,j)+1} 
$$
Supposons désormais $\k = \R$. Alors $\sqrt{\zeta(T)}$ est bien défini
parce que, sous ces hypothèses, $d_T(i,j)>0$. A
tout tableau de Young standard de forme $\la$ on associe
alors $\tilde{T} = \sqrt{\zeta(T)} T$. L'action de $s_r$ pour $1 \leq r \leq
n-1$ dans la base formée des $\tilde{T}$ pour $T$ standard
est alors donnée par 
par $s_r . \tilde{T} = \tilde{T}$ si $r$ et $r+1$ sont
sur la même colonne de $T$, par $s_r . \tilde{T} = -\tilde{T}$
si $r$ et $r+1$ sont
sur la même ligne de $T$; dans les autres cas, si $T<T_r$ on
a $d = d_T(r+1,r)>0$ et on vérifie
facilement que $\zeta(T_r) = \zeta(T) \frac{d-1}{d+1}$.
L'action de $s_r$
est alors donnée dans la base $(\tilde{T},\tilde{T}_r)$ du plan que ces éléments
engendrent
par la matrice
$$
\frac{1}{d} \left( \begin{array}{cc} -1 & \sqrt{d^2-1} \\
\sqrt{d^2-1} & 1 \end{array} \right)
$$

\paragraph{Règle de Young.}

Pour tout $n \geq 2$, on identifiera $\mathfrak{S}_{n-1}$
au sous-groupe de $\SN$ composé des permutations qui
fixent $n$. La restriction d'une représentation irréductible
$\la$ de $\SN$ au sous-groupe $\mathfrak{S}_{n-1}$ donnée
par la règle suivante, dite règle de Young :
$$
\mathrm{Res}_{\mathfrak{S}_{n-1}} \la = \bigoplus_{\mu \nearrow
\la} \mu = \bigoplus_{\mu \in
P(\la)} \mu.
$$
En particulier, cette restriction est sans multiplicités
et est composée de $\delta(\la)$ composantes
irréductibles.

\paragraph{Représentations de dimension 1.}
Les représentations de dimension 1 sont la re\-pré\-sen\-ta\-tion
triviale et la représentation donnée par la signature,
qui correspondent respectivement aux partitions
$[n]$ et $[1^n]$. On notera également $\eps$ la
représentation signature. Une propriété classique
importante pour nous est $\la' \simeq \la \otimes \eps$.

\paragraph{Représentation de réflexion.} La partition $\alpha = [n-1,1]$ correspond à une représentation
remarquable de $\SN$, appelée représentation de réflexion.
Il est classique que ses puissances alternées sont également
irréductibles, et correspondent aux représentations irréductibles
dont les diagrammes de Young sont \og en équerres \fg\ :
$$
\Lambda^p \alpha \simeq [n-p,1^p]
$$
pour $1 \leq p \leq n-1$.

\subsection{Décalages et formes bilinéaires}

Par des moyens combinatoires, on définit dans cette section
quand $\la = \la'$ une
forme bilinéaire sur l'espace vectoriel sous-jacent.

A tout tableau standard $T$ associé à une partition $\la$ de $n$, on associe la quantité suivante
$$
w(T) = \prod_{\stackrel{i<j}{c_i(T)> c_j(T)}} (-1)
= (-1)^{\# \{ i<j \ | \ c_i(T) > c_j(T) \}}
$$
Remarquons que si, pour $1 \leq r \leq n-1$, les nombres
$r$ et $r+1$ ne se trouvent ni sur
la même ligne ni sur la même colonne de $T$, alors
$w(T_r) = -w(T)$. De plus,
$$
w(T') = \prod_{\stackrel{i<j}{c_i(T')> c_j(T')}} (-1)
= \prod_{\stackrel{i<j}{l_i(T)> l_j(T)}} (-1) = (-1)^{\# \{ i<j \ | \ l_i(T) > l_j(T) \}}
$$

Repérons les cases d'un diagramme de Young associé à une partition $\la$
par les coordonnées $(i,j)$, avec $1 \leq i \leq \la'_1$,
et $1 \leq j \leq \la_i$. Il y a donc $|\la| = n$ cases.
On appelle \emph{décalage} un couple $(c_1,c_2)$
de cases $c_1 = (i_1,j_1)$ et $c_2 = (i_2,j_2)$ de $\la$,
avec $i_1 < i_2$ et $j_1 > j_2$. Le nombre de tels décalages est donc
une fonction de $\la$, notée $\nu(\la)$. A un tel décalage
on associe également le nombre positif $(d+1)/(d-1)$, avec $d
= |i_1-i_2+l_2-l_1|$, et à $\la$ le produit $\xi(\la)$ de
tous ces nombres.

Cela nous permet de démontrer le
\begin{lemme}
Si $\la$ est un diagramme de Young symétrique, pour
tout tableau de Young standard $T$ on a $w(T)w(T') = \nu(\la)$ et
$\zeta(T) \zeta(T') = \xi(\la)$.
\end{lemme}
\begin{proof}
Le produit $w(T)w(T')$
est encore \'egal au produit de $(-1)$ sur tous les couples $i<j$ tels
que $c_i(T) > c_j(T)$
ou bien, exclusivement, $l_i(T) > l_j(T)$.
De tels couples sont en bijection avec les décalages du diagramme
$\la$, de la façon suivante : au couple $(x,y)$ de cases
constituant le décalage
on associe le couple $(i,j)$ o\`u $i$ (resp. $j$) est
le minimum (resp. le maximum) de leurs contenus dans $T$ ; inversement,
la condition sur $i$ et $j$ signifie que les deux cases qui les
contiennent forment un décalage. Le deuxième énoncé se démontre
de façon similaire.
\end{proof}

On définit une forme bilinéaire non dégénérée $(\ | \ )$ sur $\la$ par
$$
(S|T) =   w(T) \delta_{S,T'}
$$
pour $S$ et
$T$ deux tableaux standards de forme $\la$, où $\delta_{i,j}$ est le symbôle de
Kronecker.
Si $S$ et $T$ sont deux tableaux standards,
$$
(S|T) = w(T) \delta_{S,T'} = w(S') \delta_{S',T} = \nu(\la) w(S) \delta_{S',T}
= \nu(\la) (T|S)
$$ 
donc $(\ | \ )$ est symétrique si $\nu(\la) =1$, et antisymétrique
si $\nu(\la) = -1$.

La forme bilinéaire $(\ | \ )$ est donc soit symétrique, soit
antisymétrique, en fonction de $\la$. Plus précisément,
quand elle est symétrique, le changement de base $U = \frac{T}{2} + T'$,
$V = \frac{T}{2} - T'$ sur chacun des plans engendrés par des couples $(T,T')$ de tableaux
standards symétriques l'un de l'autre montre que
$(\ | \ )$ est équivalente sur $\Q$ à
la forme bilinéaire symétrique standard de signature $(\frac{N}{2},\frac{N}{2})$.

Le lemme suivant simplifie le calcul de $\nu(\la)$. On rappelle
que $b(\la)$ désigne la longueur de la diagonale de $\la$.

\begin{lemme} Si $\la \vdash n$ est tel que $\la = \la'$,
le nombre de décalages de $\la$ est congru à $\frac{n - b(\la)}{2}$
modulo 2. \end{lemme}
\begin{proof}
On appelle
diagonale de $\la$ l'ensemble des cases de la forme $(i,i)$.
Il y a $b(\la)$ telles cases.

Une case $(i,j)$ est dite au dessus (resp.
au dessous) de la diagonale si $i < j$ (resp. $i>j$). On note
$\la^+$ et $\la^-$ les ensembles de cases correspondants.
Ici, $\la = \la'$ donc les ensembles $\la^+$ et $\la^-$ sont
en bijection par $\tau : (i,j) \mapsto (j,i)$.

D'autre part,
à tout décalage $c = (c_1,c_2)$ on peut associer le
décalage symétrique $\tau(c) = (\tau(c_2), \tau(c_1) )$. On en
déduit que le nombre de décalages est congru modulo 2 au
nombre de décalages fixés par $\tau$. Or, si $c = ( c_1,c_2 )$
et $\tau(c) = c$, on a $c_2 = \tau(c_1)$ et $c_1 = \tau(c_2)$. Il y a
autant de tels décalages que d'éléments de $\la^+$, soit $(n -
b(\la))/2$.
\end{proof}

Cette forme bilinéaire est reliée à l'isomorphisme
entre $\la$ et $\la \otimes \eps$. Plus précisément, introduisons
sur $\la$ le produit scalaire $<S,T> = \frac{1}{\zeta(T)} \delta_{S,T}$,
et notons $M$ l'endomorphisme de $\k$-espace vectoriel $M$ sur $\la$
par $M(T) = w(T) \zeta(T')T'$,
pour tout tableau standard $T$ de forme $\la$. Il n'est pas
difficile de vérifier que $(S|T) = <S,M(T)>$. Le produit scalaire
$<.,.>$ a la propriété que, si $\k = \R$, alors $<\tilde{S},\tilde{T}>
= \delta_{S,T}$.

\begin{lemme} \label{invfbilin}
Pour tout $s \in \SN$ et tous $x,y$ dans $\la$,
on a $M(s.x) = \eps(s) s.M(x)$,
et $(s.x|y) = \eps(s) (x|s^{-1}.y)$.
\end{lemme}
\begin{proof}
Comme ces propriétés sont linéaires en $x$ et $y$, on peut supposer,
d'une part $\k = \R$, d'autre part que $x = \tilde{S}$ et $y = \tilde{T}$
où $S$ et $T$ sont des tableaux standards. D'autre part, comme les
transpositions consécutives engendrent $\SN$, on peut supposer
$s = s_r$, pour $1 \leq r < n$. D'après les matrices du modèle
orthogonal, l'action de $s_r$ est autoadjointe pour $<\ ,\ >$,
donc orthogonale puisque $s_r^2 = 1$. Il s'ensuit que la
deuxième égalité découle de la première. D'autre
part, on a $M( \tilde{T}) = w(T) T'$. On en déduit facilement que 
$M(s_r.\tilde{T}) = -s_r.M(\tilde{T})$ : par exemple, si $T_r$
est un tableau standard et $d = d_T(r+1,r)>0$, alors d'une part
$$
M(s_r.\tilde{T}) = \frac{-1}{d} w(T) \tilde{T}' + \frac{\sqrt{d^2-1}}{d} w(T_r) \tilde{T}'_r
= \frac{-1}{d} w(T) \tilde{T}' - \frac{\sqrt{d^2-1}}{d} w(T) \tilde{T}'_r
$$
parce que $w(T_r) = -w(T)$ et, d'autre part, puisqu'alors $T' > T'_r$,
on a
$$
s_r.M(\tilde{T}) = w(T) \sqrt{d^2-1}{d} \tilde{T}'_r + w(T) \frac{1}{d} 
\tilde{T}' = -M(s_r.\tilde{T})
$$
\end{proof}

\subsection{Dimensions}

Soit $\la$ une partition de $n$, et $r \in \N$ tel que
$\la_r > 0$ et $\la_{r+1} = 0$.
La dimension du $\k$-espace vectoriel associé à $\la$,
c'est-à-dire le nombre de tableaux standards de forme
$\la$, est donné par la formule suivante (cf. \cite{FH} p. 50)
$$
\dim(\la) = \frac{n!}{l_1!\dots l_r!} \prod_{i<j}(l_i-l_j)
$$
où $l_i = \la_i +r-i$. Cette formule permet de
calculer facilement la dimension d'une représentation donnée.
Le cas des équerres se déduit immédiatement de l'isomorphisme
avec les puissances alternées de la représentation
de réflexion :
$$
\dim [n-p,1^p] = \left( \begin{array}{c} n-1 \\ p \end{array}
\right)$$
pour $0 \leq p \leq n-1$.  
Pour une approche inductive, il est souvent plus
utile d'utiliser la règle de Young, qui a pour conséquence
la formule itérative suivante :
$$
\dim([1]) = 1 \ \ \ \ \ \ \ \ \  \dim(\la)= \sum_{\mu \nearrow \la} \dim(\mu) = \sum_{\mu \in P(\la)} \dim(\mu)
$$

On sait qu'il n'y a que deux représentations de
dimension 1. On remarque ainsi
$$
\dim(\la) = 1 \Leftrightarrow \la \in \{ [n], [1^n] \} \Leftrightarrow
\delta(\la) = 1
$$
Après celles-ci, les représentations irréductibles
les plus petites sont de dimension $n-1$ :
\begin{lemme} \label{ppdim} Soit $n \geq 5$ et $\la$ une partition de $n$. Alors
\begin{enumerate}
\item $\dim(\la) > 1 \Rightarrow \dim(\la) \geq \dim(\alpha) = n-1$
\item $\dim(\la) = n-1 \Leftrightarrow \la \in \{\alpha, \alpha' \}$
\end{enumerate}
En particulier, si $\la$ est propre, 
$\dim(\la) > n-1$.
\end{lemme}
\begin{proof}
On sait que $\dim(\alpha) = \dim(\alpha') = n-1$. On procède par
récurrence sur $n$. Le cas $n=5$ se vérifie facilement par
calcul des dimensions de chacune des représentations. Si $n >5$,
soit $\la \vdash n$ avec $\dim(\la) > 1$. On a donc
$\delta(\la) > 1$, donc $P(\la)$ contient au moins
deux
éléments.

S'il existe $\mu_1 \nearrow \la$ tel que
$\dim(\mu)>1$, alors par hypothèse de récurrence on
a $\dim(\mu) \geq n-2$. Notant alors $\mu_2 \neq \mu_1$
tel que $\mu_2 \nearrow \la$ on a alors
$\dim(\la) \geq \dim(\mu_1) + 1 = n-1$. De plus, le
cas d'égalité signifie que $\mu_1 \in \{ [n-2,1],[2,1^{n-3}] \}$
et $P(\la) = \{ \mu_1,\mu_2 \}$
avec $\dim(\mu_2) = 1$, soit $\mu_2 \in \{ [n-1],[1^{n-1}] \}$.
On déduit alors facilement du fait
que $\la = \mu_1 \cup \mu_2$ qu'alors $\la \in \{\alpha,\alpha' \}$.

Au contraire, si tous les $\mu \in P(\la)$
sont de dimension 1, comme il y en a au moins deux on a
$P(\la) = \{ [n-1], [1^{n-1}] \}$
et $\la = [n] \cup [1^n] = [n-1,1^{n-2}]$, qui est une partition de 
$2n-3 > n$ pour $n \geq  4$. Ce cas est donc exclu et la
proposition démontrée.
\end{proof}

Une inégalité qui nous sera utile est la suivante :

\begin{lemme}
\label{dimdelta}
Si $|\la| \geq 6$ et $\dim \la \neq 1$, alors $\dim \la > 2 \delta(\la)$.
\end{lemme}
\begin{proof}
D'après la règle de Young, la dimension de $\la$ est la somme des
dimensions des $\mu \in P(\la)$,
qui sont au nombre de $\delta(\la)$. Pour un tel $\mu$, et
d'après le lemme \ref{ppdim}, si
$\dim(\mu) \neq 1$ alors $\dim(\mu) \geq 4 > 2$ puisque $|\mu| \geq 5$. Ainsi, on a bien
$\dim(\la) > 2 \delta(\la)$ si $\dim(\mu) \neq 1$ pour
tout $\mu \in P(\la)$. Dans le cas contraire, quitte à
échanger $\la$ en $\la'$ on peut supposer $[n-1] \in P(\la)$. Or
les seules $\la \vdash n$ telles que $[n-1] \nearrow \la$
sont $[n]$, qui est de dimension 1, et $\alpha$. On en déduit
d'une part que $\dim(\la) = \dim(\alpha) = n-1$ et d'autre part
que $\delta(\la) = 2$, et on a bien $n-1 > 2 \times 2 = 4$
puisque $n>5$.
\end{proof}

\subsection{Action de la somme des transpositions}

Une autre donnée naturellement associée à une partition
$\la$ est la valeur $\gamma(\la)$ du caractère irréductible associé sur l'une des
transpositions. Pour la calculer, on a une formule dûe à Frobenius
$$
\gamma(\la) = \frac{\dim(\la)}{n(n-1)} \sum_{i=1}^{b(\la)} \left(v_i(v_i+1) - u_i(u_i+1)
\right) = \frac{\dim(\la)}{n(n-1)} \sum_{i=1}^{b(\la)} (v_i - u_i)(v_i + u_i + 1)
$$
avec $\la_i = i + u_i$ et $\la'_i = i + v_i$ pour $1 \leq i \leq b(\la)$
(cf. \cite{FH} p. 52, ex. 4.17). En particulier,
$\gamma(\la)=0$ si $\la = \la'$.

Il nous sera utile de connaitre l'action sur chacune
des représentations de l'élément
$$
T_n = \sum_{i<j} (i \ j) \in \k \SN
$$
qui est central dans l'algèbre de groupe de $\SN$.
Comme chaque $\la$ est une représentation absolument irréductible
de $\SN$, l'élément $T_n$ agit sur elle par un scalaire. Il suffit
donc de connaître la trace de cette action, c'est-à-dire
$n(n-1)/2$ fois la valeur du caractère associé sur une transposition,
pour connaître la valeur de $T_n$. En particulier,

\begin{lemme} \label{lemactTN} L'élément $T_n \in \k \SN$ agit
sur $\la \vdash n$ par le scalaire $n(n-1)\gamma(\la)/2\dim(\la)$.
\end{lemme}

\section{Factorisations}

Dans cette section, nous associons en 5.2 à chaque partition symétrique $\la$
une algèbre de Lie $\osp(\la)$, et 
construisons les morphismes injectifs
$$
\phi_n : \g_n' \into \sl_{n-1}(\k) \times \left( \prod_{\la \in E_{n}/\sim} \sl(\la) \right) \times \left( \prod_{
\la \in F_{n}} \osp(\la) \right)
$$
mentionnés dans l'énoncé du théorème A. Par construction $\g'_n$ est une sous-algèbre de Lie semi-simple
de $\k \SN$. Or, en tant qu'algèbre de Lie,
$$
\k \SN = \bigoplus_{\la \vdash n} \gl(\la)
$$
et chaque $\la$ correspond à une représentation de $\g_n$ donc
de $\g'_n$, qui est encore irréductible parce que les
transpositions engendrent $\SN$ : il s'ensuit que tout sous-espace vectoriel
de l'espace vectoriel associé à $\la$ qui est stable par $\g'_n$ sera
stable par $\SN$. Le morphisme $\phi_n$ correspondra à une
factorisation de cette inclusion. Du simple fait que
$\g'_n$ est semi-simple, on en déduit tout d'abord qu'elle
est incluse dans la sous-algèbre de Lie dérivée de $\k \SN$, soit
la factorisation
$$
\xymatrix{
 &  {\displaystyle \bigoplus_{\la \vdash n} \gl(\la)} \\
\g'_n \ar[ur] \ar[r] & {\displaystyle \bigoplus_{\stackrel{\la \vdash n}{\dim(\la) > 1}} \sl(\la)} \ar[u] 
}
$$

\subsection{Les algèbres de Lie $\sl(\la)$ et $\sl(\la')$}

Soit $\la \vdash n$ une partition non symétrique de $n$. On note ici, pour
clarifier l'exposé, $V_{\la}$
le $\k$-espace vectoriel de base les tableaux standards de forme $\la$,
c'est-à-dire le $\k$-espace vectoriel associé à la représentation
$\rho_{\la} : \k \SN \to \End(V_{\la})$. On sait
qu'existe un isomorphisme de $\k \SN$-modules
$P_{\la} : \la' \to \la \otimes \eps$, c'est-à-dire un isomorphisme
de $\k$-espaces vectoriels entre $V_{\la'}$ et $V_{\la}$ tel que
$P_{\la}(s.x) = \eps(s)s. P_{\la}(x)$ pour tout $s \in \SN$ et
$x\in V_{\la'}$.

En tant que représentations de $\g_n$, on en déduit
que, pour toute transposition $s$ de $\SN$, donc pour
tout $s \in \g_n$, le diagramme suivant commute
$$
\xymatrix{
V_{\la'} \ar[r]^{P_{\la}} \ar[d]_{\rho_{\la'}(s)} & V_{\la} \ar[d]^{-\rho_{\la}(s)^{\#}} \\
V_{\la'} \ar[r]^{P_{\la}} & V_{\la} 
}
$$
où l'on a noté $u^{\#}$ l'adjoint de $u \in \gl(\la)$ pour un
produit scalaire $\SN$-invariant sur $\la$.
Cela signifie que $\rho_{\la}$ et $\rho_{\la'}$ sont duales l'une de
l'autre en tant que représentations de $\g_n$ donc de $\g'_n$. En particulier
$\g_{\la} \simeq \g_{\la'}$ et de plus le morphisme d'algèbres de Lie $\g'_n \to \sl(\la) \oplus \sl(\la')$ se factorise
ainsi

$$
\xymatrix{
 &  \sl(\la) \oplus \sl(\la') \\
\g'_n \ar[ur] \ar[r] & \sl(\la) \ar[u]_{\mathrm{Id} \oplus \widetilde{P_{\la}}} 
}
$$
où l'on a noté $\widetilde{P_{\la}}$ l'automorphisme de $\gl(\la)$ défini
par $\widetilde{P_{\la}}(x) = - P_{\la}^{-1} x^{\#} P_{\la}$. Il envoie naturellement
$\sl(\la)$ sur $\sl(\la')$.
On en déduit la factorisation suivante
$$
\xymatrix{
 &  {\displaystyle \bigoplus_{\stackrel{\la \vdash n}{\dim(\la) > 1}} \sl(\la)} \\
\g'_n \ar[ur] \ar[r] & {\displaystyle \left( \bigoplus_{\stackrel{\dim(\la)>1}{\la < \la'}} \sl(\la) \right)
\oplus \left( \bigoplus_{\stackrel{\dim(\la)>1}{\la = \la'}} \sl(\la) \right)
} \ar[u] 
}
$$

\subsection{L'algèbre de Lie $\osp(\la)$}

Soit $\la$ une partition symétrique de $n$. On a construit
en 4.2 une forme bilinéaire $(\ | \ )$ non dégénérée, orthogonale
ou symplectique suivant $\la$, sur l'espace vectoriel
sous-jacent. Elle nous permet de définir l'algèbre
de Lie suivante.

\begin{defi} On note $\osp(\la)$ la sous-algèbre de
Lie de $\gl(\la)$ composée des $m \in \gl(\la)$
tels que $(m.x|y) + (x|m.y) = 0$
pour tous $x,y$.
\end{defi}

D'après la section 4.2 cette algèbre de Lie simple, orthogonale
ou symplectique, est dans tous les cas de rang semi-simple
$\dim(\la)/2$.
L'intérêt pour nous de cette algèbre de Lie vient du fait suivant.

\begin{lemme} \label{glaosp} Pour tout $\la \vdash n$ symétrique,
$\g_{\la} \subset \osp(\la)$.
\end{lemme}
\begin{proof} 
Comme $\g_n$ est engendrée par les transpositions,
il suffit de montrer que l'image de chacune d'entre elles
appartient à $\osp(\la)$, c'est-à-dire
que $(s.x|y) = -(x|s.y)$ pour tous $x,y$ dans $\la$.
C'est une conséquence immédiate du lemme \ref{invfbilin} puisque,
lorsque $s$ est une réflexion, $s^{-1} = s$ et $\eps(s) = -1$.
\end{proof}

On en déduit la factorisation

$$
\xymatrix{
 &  {\displaystyle \left( \bigoplus_{\stackrel{\dim(\la)>1}{\la < \la'}} \sl(\la) \right)
\oplus \left( \bigoplus_{\stackrel{\dim(\la)>1}{\la = \la'}} \sl(\la) \right)
} 
 \\
\g'_n \ar[ur] \ar[r] & {\displaystyle \left( \bigoplus_{\stackrel{\dim(\la)>1}{\la < \la'}} \sl(\la) \right)
\oplus \left( \bigoplus_{\stackrel{\dim(\la)>1}{\la = \la'}} \osp(\la) \right)
} \ar[u] 
}
$$

\subsection{Equerres}

On note $\alpha_r = [n-r,1^r]$ pour tout $1 \leq r \leq n-1$, 
ainsi $\alpha = \alpha_1 = [n-1,1]$. On a rappelé en section 4.1 que,
en tant que représentations de $\mathfrak{S}_n$,
$\Lambda^r \alpha = \alpha_r$, où $\Lambda^r V$ désigne la
puissance extérieure $r$-ième\ de $V$. Un phénomène remarquable
est que les représentations de $\g'_n$ correspondantes sont également
isomorphes.
Ce résultat découle immédiatement de la classification
des systèmes KZ irréductibles sous l'action du groupe symétrique
(cf. \cite{THESE,KZ}). Nous en donnons ici une démonstration
élémentaire :

\begin{lemme} \label{isoequerres} En tant que représentation de $\g'_n$, $\alpha_r
= \Lambda^r \alpha$.
\end{lemme}
\begin{proof}
Soit $\beta$ la représentation naturelle de $\SN$, définie sur
une base $v_1,\dots,v_n$ par $s. v_i = v_{s(i)}$ pour tout
$s \in \SN$. On a $\beta = [n] \oplus \alpha$ en tant que représentation
de $\SN$ et de $\g_n$. Alors $\Lambda^r \beta = \Lambda^r \alpha
\oplus \Lambda^{r-1} \alpha$. Comme $\Lambda^{n-1} \beta = \Lambda^{
n-1} \alpha \oplus \Lambda^{n-2} \alpha$, il suffit de montrer
que l'action de $\g'_n$ sur $\Lambda^r \beta$ est déduite de celle
de $\SN$ par la composée $\g'_n \to \g_n \to \k \SN$ pour en déduire
$\alpha_r = \Lambda^r \alpha$ en tant que représentation de $\g'_n$.
Fixons $r \in [1,n]$, et notons $\underline{v} = v_{i_1}\wedge \dots
\wedge v_{i_r}$ pour $i_1 < i_2< \dots < i_r$, ainsi que $\tau_{ij} = (i \ j)
\in \g_n$ et $\sigma_{ij} = (i \ j) \in \SN$.
\begin{itemize}
\item Si $i_1 \geq 3$, on a $\tau_{12} . \underline{v}
= r \underline{v}$ et $\sigma_{12} . \underline{v} = \underline{v}$.
\item Si $i_1 = 1$ et $i_2 \geq 3$, ou $i_1 = 2$ et
$i_2 \geq 3$, on a $\tau_{12} . \underline{v} = \sigma_{12}. \underline{v}
+ (r-1) \underline{v}$.
\item Si $i_1 = 1$ et $i_2 = 2$, on a $\tau_{12}. \underline{v}
= (r-2) \underline{v}$ et $\sigma_{12}.\underline{v} = - \underline{v}$
\end{itemize}
et dans tous les cas
on a donc $\tau_{12} . \underline{v} = \sigma_{12}. \underline{v}
+ (r-1) \underline{v}$. On en déduit que l'action de $\tau_{12}$
est celle de $\sigma_{12} + (r-1) \Id$. Comme l'action de
$\g_n$ est $\SN$-équivariante pour l'action par conjugaison de
$\SN$ sur $\g_n \subset \k \SN$, l'action de $\tau_{ij}$
est celle de $\sigma_{ij} + (r-1) \Id$, et l'action de $\g'_n$
sur $\Lambda^r \beta$ est bien déduite de l'action de $\SN$
par $\g'_n \subset \k \SN$, d'où $\alpha_r = \Lambda^r \alpha$
en tant que représentation de $\g'_n$.
\end{proof}

On en déduit que les morphismes $\g'_n \to \sl(\alpha_r)$
se factorisent en
$$
\xymatrix{
 & \sl(\alpha_r) \\
\g'_n \ar[ur]^{\rho_{\alpha_r}} \ar[r]_{\rho_{\alpha}} &
\sl(\alpha) \ar[u]^{\Delta_r} }
$$
où $\Delta_r$ est donné par
$$
\Delta_r(x) = \sum_{i=1}^r 1 \wedge \dots \wedge \stackrel{(i)}{x} \wedge \dots \wedge 1
$$
On en déduit la factorisation finale
$$
\xymatrix{
 &  {\displaystyle \left( \bigoplus_{\stackrel{\dim(\la)>1}{\la < \la'}} \sl(\la) \right)
\oplus \left( \bigoplus_{\stackrel{\dim(\la)>1}{\la = \la'}} \osp(\la) \right)
} 
 \\
\g'_n \ar[ur] \ar[r] & \sl_{n-1}(\k) \times \left( \prod_{\la \in E_{n}/\sim} \sl(\la) \right) \times \left( \prod_{
\la \in F_{n}} \osp(\la) \right) \ar[u] 
}
$$

\section{Démonstration du théorème A}

\subsection{Réduction à $\k = \C$}

Si $\k$ est un corps de caractéristique 0, notons $\g_n(\k) \subset \k
\SN$ l'algèbre de Lie des transpositions définie sur $\k$, et
$$
\mathcal{L}_n(\k) = 
\sl_{n-1}(\k) \times \left( \prod_{\la \in E_{n}/\sim} \sl(\la) \right) \times \left( \prod_{
\la \in F_{n}} \osp(\la) \right) \subset \k \SN
$$
l'algèbre de Lie de dimension finie construite en section 5.3. On a construit,
pour tout corps $\k$ de caractéristique 0, un morphisme
injectif $\phi_n^{\k} : \g_n(\k) \to \mathcal{L}_n(\k)$.
Par construction, la dimension de $\mathcal{L}_n(\k)$ sur $\k$
est indépendante de $\k$, et celle de $\g_n(\k)$ également
parce que $\g_n(\k)  = \g_n(\Q) \otimes_{\Q} \k$. Il suffit
donc de montrer l'égalité des dimensions pour $\k = \C$, c'est-à-dire
de montrer que $\phi_n^{\C}$ est surjectif, pour conclure la
preuve du théorème.

A partir de maintenant, on supposera donc $\k = \C$. Toutes les
algèbres de Lie considérées sont donc semi-simples complexes,
donc déterminées à isomorphisme près par leur type dans
la classification de Cartan. On note à nouveau
$\g_n = \g_n(\C)$, et on notera $\mathcal{L}_n = \mathcal{L}_n(\C)$.

\subsection{Généralités sur les algèbres de Lie semi-simples}

Nous utiliserons la rigidité dimensionnelle de la classification des
algèbres de Lie simples complexes, notamment
les deux lemmes suivants, qui se déduisent facilement
de cette classification et de la formule des
caractères de Weyl par examen des cas (cf. par exemple les calculs
de \cite{FH}, ex. 24.52). On utilise,
pour chaque algèbre de Lie simple
de rang $n$,
la numérotation des poids fondamentaux $\varpi_1,\dots, \varpi_n$
de \cite{FH}, et on accepte les redondances dans la classification
des algèbres de Lie simples ($A_3 = D_3$, $B_2 = C_2$, etc.).
Pour toute algèbre semi-simple complexe $\g$, on notera
$\rg(\g)$ son rang.

\begin{lemme} \label{caractA} Si $\g$ est une algèbre de Lie simple
complexe qui admet une représentation $V$ telle
que $\dim(V) < 2 \rg(\g)$, alors $\g \simeq \sl(V)$.
\end{lemme}

\begin{lemme} \label{caractBCD} Les paires $(\g,V)$ avec $\g$ une algèbre de Lie
simple complexe de rang $n$ et $V$ une représentation irréductible
de $\g$ de dimension $N$ telles que $2n \leq N < 4n$ sont,
pour tout $n \geq 2$, $(B_n,\varpi_1)$ pour $N = 2n+1$, $(C_n,\varpi_1)$
pour $N = 2n$,
et $(D_n,\varpi_1)$ pour tout $n \geq 3$ ($N = 2n$) plus, pour $n \leq 6$,
les couples exceptionnels suivants :
\begin{itemize}
\item[Rang 1 :] $(A_1,\varpi_1)$ pour $N = 2$, $(A_1,2\varpi_1)$ pour
$N = 3$.
\item[Rang 2 :] $(G_2,\varpi_1)$ pour $N = 7$, $(B_2,\varpi_2)$ pour
$N = 4$, $(C_2,\varpi_1)$ pour $N = 5$.
\item[Rang 3 :] $(B_3,\varpi_3)$ pour $N = 8$, $(A_3,\varpi_2)$ pour
$N = 6$.
\item[Rang 4 :] $(D_4,\varpi_3)$ et $(D_4,\varpi_4)$ pour $N = 8$, $(A_4,
\varpi_2)$ et $(A_4,\varpi_3)$ pour
$N = 10$.
\item[Rang 5 :] $(A_5,\varpi_2)$ et $(A_5,\varpi_4)$ pour $N = 15$, $(D_5,\varpi_4)$
et $(D_5,\varpi_5)$pour
$N = 16$.
\item[Rang 6 :] $(A_6,\varpi_2)$ et $(A_6,\varpi_5)$ pour
$N = 21$.
\end{itemize}
\end{lemme}

Pour se ramener au cas des algèbres de Lie simples
et utiliser les deux lemmes précédents, nous invoquerons le
lemme suivant.

\begin{lemme} \label{lemredsimp} Soit $V$ un $\C$-espace vectoriel de dimension finie,
$\h \subset \g$ deux sous-algèbres de Lie semi-simples de $\sl(V)$.
On note $\h = \bigoplus_{j \in J} \h_j$ la décomposition de $\h$
en idéaux simples. On suppose que les propriétés suivantes
sont vérifiées
\begin{enumerate}
\item $V$ est irréductible sous l'action de $\g$.
\item Pour tout $j \in J$, la restriction de $V$ à $\h_j$ admet
une composante irréductible de multiplicité 1.
\item On a $\rg(\g) < 2 \rg(\h)$.
\end{enumerate}
Alors $\g$ est une algèbre de Lie simple.
\end{lemme}
\begin{proof}
On note $\g = \g^1 \oplus \dots \g^g$ la décomposition de $\g$ en idéaux
simples, et on fixe un $i \in [1,g]$. On note $p:\g \onto \g^i$
la projection sur $\g^i$ parallèlement à
$\g^{(i)} = \bigoplus_{r \neq i} \g^r$. Soit $q : \h \to \g^i$ la
restriction de $p$ à $\h$. On montre par l'absurde que
$q$ est injective.

En effet, $\Ker q$ serait sinon un idéal non nul de $\h$, donc une
somme non vide $\bigoplus_{j \in K} \h_j$ de certains de ses
idéaux simples, avec $K \subset J$ et $K \neq \emptyset$. On choisit
un tel $j \in K$. Comme $V$ est irréductible et fidèle pour
l'action de $\g$, cette représentation est isomorphe
à une représentation de la forme $U_i \otimes U_{(i)}$ avec
$U_i$ une représentation fidèle de $\g^i$ et
$U_{(i)}$ une représentation de $\g^{(i)}$. Comme
$\h_j \subset \Ker q$, cet idéal agit trivialement sur
$U_i$, donc chacune des composantes irréductibles de la restriction
à $\h_j$ de $V$ intervient avec multiplicité au moins $\dim U_i$.
Cela contredit l'hypothèse 2) parce que $\g^i$ est non commutative donc
sa représentation fidèle $U_i$ est de dimension au moins 2.

Ainsi, $\h$ s'injecte dans $\g^i$ pour tout $i \in [1,g]$,
donc $\rg (\g) \geq g \rg (\h)$. Alors $g \leq \rg (\g) / \rg (\h) < 2$
donc $g = 1$ et $\g$ est une algèbre de Lie simple.
\end{proof}

\paragraph{Remarque.} L'hypothèse 3) du lemme précédent
est notamment vérifiée lorsque $\rg (\h) \geq (\dim V)/2$. Elle est
également vérifiée lorsque
$\g \subset \so(V)$ ou $\g \subset \sp(V)$ et $\rg (\h) > (\dim V)/4$.

\subsection{Diagrammes remarquables}

\begin{figure}
\begin{center}
\includegraphics{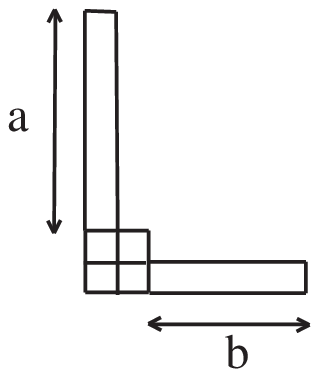}
\end{center}
\caption{Les diagrammes $\D(a,b)=[a+2,2,1^b]$}
\label{figdab}
\end{figure}

Pour tous $a,b \geq 0$, on note $\D(a,b)= [a+2,2,1^b]$. Il s'agit d'une partition de
$n = a+b+4$, dont le diagramme est de la forme décrite dans la figure \ref{figdab}.

D'après la règle de Young, pour $\la \vdash n$, l'ensemble $P(\la)$ contient une équerre,
c'est-à-dire un $\mu \vdash n-1$ tel que $b(\mu)=1$, si et seulement si $\la$
est elle-même une équerre ou bien si $\la$ est justement de cette forme. Plus précisément,
si $\la = \D(a,b)$ avec $a,b \geq 1$, on a
$$
P(\D(a,b)) = \{ \D(a-1,b) , \D(a,b-1) , [a+2,1^{b+1}] \}
$$
Remarquons également $\D(a,b)' =\D(b,a)$ et en particulier
$\dim \D(a,b) = \dim \D(b,a)$. La formule des dimensions permet d'obtenir facilement
$$
\dim \D(a,b) = \frac{b+1}{a+2} \cnp{n-2}{a} n
$$
où $\D(a,b) \vdash n$, c'est-à-dire $n = a+b+4$. En particulier, pour tout
$a \geq 1$,
$$
\dim \D(a,a) = 2(a+1) \cnp{2a+2}{a}\ \ \mbox{ et } \ \ \dim \D(a-1,a) = \cnp{2a+1}{a-1} (2a+3).
$$
et enfin, pour tout $m \geq 2$,
$$
\dim \D(m,m-2) = \frac{2(m-1)(m+1)}{m+2} \cnp{2m}{m}
$$

Les résultats techniques suivants nous seront utiles. Nous les regroupons
sous la forme d'un lemme, dont la démonstration est élémentaire et laissée au lecteur.

\begin{lemme} \label{leminegdab}
\begin{enumerate}
\item Pour tout $m \geq 3$, on a $\dim \D(m,m-2) > \cnp{2m+1}{m}$.
\item Pour tout $a \geq 2$, on a $\dim \D(a-1,a) > \frac{1}{4} \dim \D(a,a)$.
\item Pour tous $a > b+1 \geq 3$, on a $\dim \D(a,b) > 3 \cnp{a+b+2}{b+1}$
\end{enumerate}
\end{lemme}

\subsection{Résultats précédents et hypothèse de récurrence}

Les résultats suivants sur l'algèbre de Lie $\g'_n$ ont été démontré dans des travaux antérieurs

\begin{prop} \label{casTL} (\cite{QUOTDEF}, théorème 3.) Pour tout $1 \leq r \leq n$, si $\la = [n-r,r]$,
alors $\g_{\la} = \sl(\la)$.
\end{prop}

\begin{prop} (\cite{QUOTDEF}, section 8.1) \label{petitscas}
Pour tout $n \leq 6$, le morphisme $\phi_n$ est surjectif.
\end{prop}

Nous allons démontrer le théorème $A$ par récurrence sur $n$
à partir de cette dernière proposition, prenant pour
hypothèse de récurrence
$$
\mathrm{(HR_{n})}\ \ \ \ \   \phi_n \ \mathrm{est\ surjectif.}
$$

On peut donc supposer $\HRm$ vrai pour un certain $n \geq 7$.
La démonstration que dans ce cas $\HRm \Rightarrow \HRn$
fait l'objet des sections suivantes. 

\subsection{Détermination des $\g_{\la}$}

On suppose désormais $\HRm$, avec $n \geq 7$. On va dans
un premier temps en déduire $\g_{\la}$ pour
$\la \vdash n$.

On note $\h_{\la}$ l'image de $\g'_{n-1} \subset \g'_n$ dans $\g_{\la}$.
Cette image se détermine à partir du diagramme commutatif
suivant
$$
\xymatrix{
\g'_n \ar[rr] & & \sl(\la) \\
\g'_{n-1} \ar[u] \ar[rr] \ar[dr]_{\phi_{n-1}} & & \bigoplus_{\mu \nearrow \la} \sl(\mu) \ar[u] \\
 & \mathcal{L}_{n-1} \ar[ur]
}
$$
Par hypothèse de récurrence, $\phi_{n-1}$ est un isomorphisme,
donc $\h_{\la}$ est l'image de $\mathcal{L}_{n-1}$
dans les endomorphismes de sa représentation
$\bigoplus_{\mu \nearrow \la} \mu$ et se détermine aisément en fonction
de $\la$. D'autre part, tout facteur simple de $\mathcal{L}_{n-1}$
associé à un $\mu \vdash n-1$ est soit de rang $\dim(\mu)-1$,
soit de rang $\dim(\mu)/2$. A partir de la règle de Young
et de la combinatoire des diagrammes, on en déduit
les résultats suivants.

\subsubsection{Les cas particuliers $\la = \D(a,b)$}

Pour $a,b \geq 0$ et $\la = \D(a,b)$, on note $\dd_{a,b} = \h_{\la}$ et $N_{a,b} = \dim(\la)$. Si l'on convient
$N_{a,b} = 0$ pour $a<0$ ou $b <0$, la règle de Young implique
$$
N_{a,b} = N_{a-1,b} + N_{a,b-1} + \cnp{a+b+2}{b+1}
$$
En particulier, les nombres $N_{a,b}$ varient de façon croissante en $a$ et en $b$.

\begin{lemme} \label{casdab} Soient $a,b >0$ tels que $n=a+b+4 \geq 7$. Alors, sous $\HRm$, on a
\begin{enumerate}
\item Si $a=b$ alors $\rg(\dd_{a,a}) > \frac{1}{4} N_{a,a}$.
\item Si $a \neq b$, alors $\rg(\dd_{a,b}) > \frac{1}{2} N_{a,b}$.
\end{enumerate}
\end{lemme}
\begin{proof}
Si $b=0$, on a $a \geq 3$ puisque $n \geq 7$. D'après la proposition \ref{casTL},
comme $P(\D(a,0)) = \{ D(a-1,0) , [a+2,1] \}$ on a, d'une part $N_{a,a} = N_{a-1,0} + a+2$, et d'autre
part $\dd_{a,0} \simeq \sl(\D(a-1,0))\times
\sl_{a+2}(\C)$. On en déduit $\rg \dd_{a,0} = N_{a-1,0} +a = N_{a,0} - 2$. Alors $N_{a,0} - 2 > > \frac{1}{2}N_{a,0}$
si et seulement si $N_{a,0} > 4$, ce qui est vrai car $N_{a,0} \geq N_{3,0} = 14$.
Puisque $\D(a,b)' = \D(b,a)$ on
peut donc supposer $a \geq b > 0$.

Supposons d'abord $a=b > 0$. Alors $P(\D(a,b)) = \{ \D(a-1,a) + \D(a,a-1) + [a+2,1^{a+1}] \}$.
D'après $\HRm$, on alors
$$
\rg \dd_{a,a} = N_{a-1,a} - 1 + n-4 = N_{a-1,a} + 2a -1 > N_{a-1,a} > \frac{1}{4} N_{a,a}$$
d'après la partie 2 du lemme \ref{leminegdab}, puisque $n = 2a+4 \geq 7$ implique $a \geq 2$. Cela montre 1).

Si $a > b > 0$, pour montrer 2) il y a deux cas à considérer.
\begin{enumerate}
\item Si $a=b+1$, alors $\la = \D(a,b) = \D(b+1,b)$ avec
$b \geq 1$, et $P(\la) = \{ \D(b,b) + \D(b+1,b-1) + [b+3,1^{b+1}] \}$. D'après $\HRm$,
$$
\rg(\dd_{b+1,b})  =  \frac{N_{b,b}}{2} + N_{b+1,b-1} + n-4 
= \frac{ N_{b+1,b}}{2} + \frac{N_{b+1,b-1}}{2} - \frac{\dim [b+3,1^{b+1}]}{2} + n-4 
$$
Posant $m = b+1$ on déduit de la partie 1) du lemme \ref{leminegdab} que, si $b \geq 2$, alors
$$
\rg(\dd_{b+1,b}) > \frac{N_{b+1,b}}{2} + n-4 > \frac{N_{b+1,b}}{2}
$$
puisque $\dim [b+3,1^{b+1}] = \cnp{2b+3}{b+1}$. On traite séparément le cas $b=1$, soit
$a =2$. Alors $N_{1,1} = 16$ et $N_{2,0} = 9$. On en déduit
$\rg (\dd_{2,1}) = \frac{N_{1,1}}{2} + (N_{2,0}-1) + 4 = 20 > N_{2,1}/2$
car $N_{2,1} = 35$.

\item Si $a > b+1 \geq 2$, alors d'après $\HRm$ on a $\rg(\g_{\D(a-1,b)}) = N_{a-1,b} - 1$
et $\rg(\g_{\D(a,b-1)}) = N_{a,b-1} - 1$. On en déduit
$$
\rg \dd_{a,b} = N_{a,b} - \dim [a+2,1^{b+1}] - 2 + n-3.
$$
Si $b+1 \geq 3$, utilisant la partie 3) du lemme \ref{leminegdab} on obtient
$$
N_{a,b} - \cnp{a+b+2}{b+1} + n-5 > N_{a,b} - \cnp{a+b+2}{b+1} > \frac{2}{3} N_{a,b} > \frac{1}{2} N_{a,b}
$$
On traite séparément le cas $b=2$. Dans ce cas, on a
$\rg (\dd_{a,2}) = \rg (\g_{\D(a-1,2)}) + \rg (\g_{\D(a,1)}) + a+1$, soit
$$
\rg \dd_{a,2} = N_{a-1,2} + N_{a,1} + (a+2) - 3 \geqslant N_{a,2} -3
$$
car $\dim [a+2,1] = a+2$. Or $N_{a,2} -3 > N_{a,2}/2$ équivaut à $N_{a,2} > 6$,
ce qui est vrai car $a \geq 4$ donc $N_{a,2} \geq N_{1,1} = 16>6$.
\end{enumerate}

\end{proof}

\subsubsection{Cas général}

On rappelle qu'une partition est dite propre si ce n'est pas une équerre, et que l'ensemble
des partitions propres de $n$ est subdivisée entre l'ensemble $F_n$ des partitions propres
symétriques et l'ensemble $E_n$ de celles qui ne le sont pas.

Pour une partition $\la$ donnée, on notera
$$
P_E(\la) = E_n \cap P(\la) \mbox{    et    } P_F(\la) = F_n \cap P(\la).
$$
Remarquons enfin que, si $\la \in E_n \cup F_n$, alors $P(\la)$ ne contient aucune représentation
de dimension 1. De plus, dire dans ce cas que $P(\la)$ contient une
équerre équivaut à dire que $\la$ est de la forme $\D(a,b)$ pour
certains $a,b \geq 0$.

\begin{prop} \label{estimrk} Soit $n \geq 7$. Si $\HRm$ est vraie, on a
$\la \in E_n \Rightarrow \rg (\g_{\la}) > \frac{\dim (\la)}{2}$ et
$\la \in F_n \Rightarrow \rg (\g_{\la}) > \frac{\dim (\la)}{4}$.

\end{prop}
\begin{proof}
Si $P(\la)$ contient une équerre, alors $\la$ est de la forme $\D(a,b)$ pour certains
$a,b \geq 0$. On en déduit la conclusion d'après la proposition \ref{casTL} si $a$ ou $b$ est nul, et
d'après le lemme \ref{casdab} sinon. On peut donc supposer que $P(\la) \subset E_n \cup F_n$.

On a $\rg (\g_{\la}) \geq \rg (\h_{\la})$ et, d'après $\HRm$, on a 
$$
\h_{\la} \simeq \left( \bigoplus_{\mu \in P_E(\la)/\sim} \sl(\mu) \right) \oplus \left( \bigoplus_{\mu \in P_F(\la)} \osp(\mu) \right)
$$

Supposons d'abord $\la \in E_n$. D'après le lemme \ref{multilemme}-3), pour $\mu \in P_E(\la)$
on a $\mu' \not\in P_E(\la)$ donc $P_E(\la)/\sim$ s'identifie à $P_E(\la)$. On déduit alors de la
décomposition précédente le rang de $\h_{\la}$ en fonction de $\delta(\la) = \# P(\la)$ et de $\dim(\la)$,
qui vaut $\sum_{\mu \in P(\la)} \dim (\mu)$
d'après la règle de Young. En effet, d'après le lemme \ref{multilemme}-4), on a $\# P_F(\la) \leq 1$. Comme
$\rg(\sl(\mu)) = \dim(\mu)-1$ et $\rg (\osp(\mu)) = \frac{\dim(\mu)}{2}$, on obtient immédiatement
\begin{enumerate}
\item Si $P_F(\la) = \emptyset$, alors $P_E(\la) = P(\la)$ et $\rg (\h_{\la}) = \dim(\la) - \delta(\la)$.
\item Si $P_F(\la) = \{ \mu_0 \}$, alors $\rg (\h_{\la}) = \dim(\la) - \frac{\dim (\mu_0)}{2} - \delta(\la) + 1$.
\end{enumerate}
Dans chacun de ces deux cas, vérifions $\rg (\h_{\la}) > \frac{\dim(\la)}{2}$. Dans le premier cas cela signifie
$\dim(\la) - \delta(\la) > \frac{\dim(\la)}{2} \Leftrightarrow \dim(\la) > 2 \delta(\la)$, ce qui est vrai d'après le lemme
\ref{dimdelta}, puisque $\la \in E_n$ implique $\dim(\la) \neq 1$. Dans le deuxième, $\rg (\h_{\la}) > \frac{\dim(\la)}{2}$
signifie $\dim(\la) - \dim(\mu_0) > 2(\delta(\la) - 1)$. Or, comme $\la$ n'est pas une équerre, tout $\mu \in P(\la)$
vérifie $\dim(\mu)>1$, et donc $\dim(\mu) \geq n-2 > 4$ d'après le lemme \ref{ppdim}. Or
$$
\dim(\la) - \dim(\mu_0) = \sum_{\mu \in P_E(\la)} \dim(\mu) > 4 \# P_E(\la) = 4(\delta(\la) - 1)$$
d'où $\dim(\la) - \dim(\mu_0) >4(\delta(\la) - 1) \geq 2(\delta(\la)-1)$, ce qui conclut.

\medskip

Supposons ensuite $\la \in F_n$. D'après le lemme \ref{lemmesym} on a $\mu' \in P(\la)$
dès que $\mu \in P(\la)$, donc $\# (P_E(\la)/\sim) = \#P_E(\la)/2$ et
$$
\sum_{\mu \in P_E(\la)/\sim} \dim(\mu) = \frac{1}{2} \sum_{\mu \in P_E(\la)} \dim(\mu)
= \frac{1}{2} \dim(\la) - \frac{1}{2} \sum_{\mu \in P_F(\la)} \dim(\mu)
$$
Puisque $\# P_F(\la) \leq 1$, on a comme précédemment deux cas. On obtient immédiatement
\begin{enumerate}
\item Si $P_F(\la) = \emptyset$, alors $\# (P_E(\la)/\sim) = \delta(\la)/2$ et
$\rg (\h_{\la}) = (\dim(\la) - \delta(\la))/2$.
\item Si $P_F(\la) = \{ \mu_0 \}$, alors $\# (P_E(\la)/\sim )= (\delta(\la)-1)/2$ et
$\rg (\h_{\la}) = (\dim(\la)-\delta(\la)+1)/2$.
\end{enumerate}
Dans chacun de ces deux cas, vérifions $\rg (\h_{\la}) > \frac{\dim(\la)}{4}$. Dans le premier
cas cela signifie $\dim(\la) > 2 \delta(\la)$, ce qui est vrai d'après le lemme \ref{dimdelta}.
Dans le deuxième cas cela signifie $\dim(\la) > 2 \delta(\la)-2$, ce qui est {\it a fortiori} vrai
pour la même raison et conclut la preuve.

\end{proof}

Pour en déduire les $\g_{\la}$ nous aurons encore besoin du lemme
suivant. Remarquons d'abord que, d'après la règle de Young, si $\mu
\subset \la$ alors $\g_{\mu}$ s'injecte dans $\g_{\la}$ comme
sous-algèbre de Lie.

\begin{lemme} \label{minork} Pour $n \geq 7$ et $\la \vdash n$, on a $\rg (\g_{\la}) \geq 4$.
Si $\la \in F_n$, on a de plus $\rg (\g_{\la}) \geq 8$.
\end{lemme}
\begin{proof}
Si $\la$ est une équerre, alors $\rg (\g_{\la}) = n-2 \geq 5$ d'après la
proposition \ref{casTL}. Si $\la \in E_n$ et $n \geq 7$, on ne
peut avoir $\la_1 < 3$ et $\la'_1 <3$, car sinon $n = \la_1 + \la_2 \leq 4$.
Comme $\g_{\la} \simeq \g_{\la}'$ on peut donc supposer $\la_1 \geq 3$.
Alors $\la_2 \geq 2$ car sinon $\la$ serait une équerre. Ainsi
$[3,2] \subset \la$ et $\g_{\la}$ contient $\g_{[3,2]}$ qui est de
rang 4 par surjectivité de $\phi_5$ (proposition \ref{petitscas}) donc
$\rg (\g_{\la}) \geq 4$. Si enfin 
$\la \in F_n$, on a donc également $\la'_1 = \la_1 \geq 3$, et
$\la_2 \geq 2$ puisque $\la$ n'est pas une équerre, ainsi
$\g_{\la}$ contient $\g_{[3,2,1]}$ qui est de rang 8 par surjectivité
de $\phi_6$ d'où
$\rg (\g_{\la}) \geq 8$.
\end{proof}

\begin{prop} \label{detgl} Soit $n \geq 7$. Si $\HRm$ est
vraie, alors $\la \in E_n \Rightarrow \g_{\la} = \sl(\la)$
et $\la \in F_n \Rightarrow \g_{\la} = \osp(\la)$.
\end{prop}
\begin{proof}
Les idéaux simples de l'algèbre de Lie semi-simple
$\h_{\la}$ sont de la forme $\sl(\mu)$, $\osp(\mu)$ pour
$\mu \nearrow \la$ et $\mu \in E_{n-1} \cup F_{n-1}$, ou bien
isomorphe à $\sl(\alpha) \simeq \sl_{n-2}(\C)$.

Dans les deux premiers cas, d'après la règle de Young,
la restriction de $\la$ à cet idéal simple comporte
toujours la représentation standard de $\sl(\mu) $ ou $\osp(\mu)$
avec multiplicité 1.

Dans le deuxième cas, c'est-à-dire si $\la$ est de la forme $\D(a,b)$,
la restriction à l'idéal simple $\sl(\alpha)$ est somme de la
puissance alternée $r$-ième de sa représentation standard,
pour $[n-1-r,1^r] = [a+2,1^{b+1}] \nearrow \D(a,b)$, et de
copies de la représentation triviale. Cette puissance
alternée étant non triviale, puisque $r = b+1 \geq 1$, elle intervient
bien avec multiplicité 1.

Le lemme \ref{lemredsimp} et la remarque qui le suit permettent
de déduire de la proposition \ref{estimrk} que $\g_{\la}$
est une algèbre de Lie simple. Si $\la \in E_n$, le fait
que $\rg( \g_{\la}) > \dim (\la)/2$ permet de déduire du lemme \ref{caractA}
que $\g_{\la} = \sl(\la)$. Si $\la \in F_n$, on a
$\rg (\g_{\la}) \geq 8$ d'après le lemme \ref{minork}. Comme
$\rg (\g_{\la}) > \dim (\la)/4$ d'après la proposition \ref{estimrk},
on déduit du lemme \ref{caractBCD} que $\g_{\la}$
est l'ensemble des éléments de $\sl(\la)$ qui laissent invariant
un certain sous-espace, de dimension 1, de $\la^* \otimes \la^*$.
Or cet espace contient les $\osp(\la)$-invariants de $\la^* \otimes
\la^*$ puisque $\g_{\la} \subset \osp(\la)$, qui forment également
un sous-espace de dimension 1. On en déduit $\g_{\la} = \osp(\la)$.
\end{proof}

\subsection{Idéaux simples de $\g'_n$}

Soit $r \in \N$, et $\la \vdash r$. On note
$\g^{\la}$ l'idéal de $\g'_r$ associé
à la représentation $\rho_{\la} : \g'_n \to \g_{\la}$,
c'est-à-dire l'orthogonal du noyau de $\rho_{\la}$ vis-à-vis
de la forme de Killing de l'algèbre de Lie semi-simple
$\g'_r$. Par définition, $\rho_{\la}$ induit un isomorphisme
d'algèbres de Lie entre $\g^{\la}$ et $\g_{\la}$.

Les factorisations de la section 5 impliquent que $\g^{\la} = \g^{\la'}$,
et que $\g^{\la} = \g^{\alpha} \simeq \sl_{n-1}(\C)$ si $\la$ est une équerre. En
tant qu'algèbre de Lie semi-simple, $\g'_r$ est produit de
ses ideaux simples.

Considérons la propriété 
$$
\mathrm{(IS_r)}\ \ \ \ \ \ \ \mbox{ Si } \la , \mu \vdash r \mbox{ sont tels que }
\{ \la, \mu \} \cap (E_r \cup F_r) \neq \emptyset,
\mbox{ alors } \g^{\la} = \g^{\mu} \Leftrightarrow \mu \in \{ \la, \la' \}$$

Le lemme suivant est immédiat.
\begin{lemme} \label{surjcomb} Si $\phi_r$ est surjectif, les idéaux simples de
$\g'_r$ sont $\g^{\alpha}$ et les $\g^{\la}$ pour $\la \in E_r \cup F_r$,
et $\mathrm{(IS_r)}$ est vérifiée.
\end{lemme}

Il admet une réciproque partielle

\begin{lemme} \label{lemconclu} Sous $\HRm$, si $\mathrm{(IS_n)}$ est vérifiée, alors
$\phi_n$ est surjectif.
\end{lemme}
\begin{proof} Comme $\g'_n$ est semi-simple, elle est produit de
ses idéaux simples. D'après la proposition \ref{detgl}, pour tout $\la \in
E_n \cup F_n$, $\g^{\la} \simeq \g_{\la}$ est simple. L'hypothèse
$\mathrm{(IS_n)}$ dit alors que l'on ne peut avoir $\g^{\la} = \g^{\mu}$
que si $\la$ et $\mu$ sont des équerres, ou bien si $\mu \in \{ \la, \la' \}$.
Ainsi,
$$
\g'_n \supset \g^{\alpha} \times \prod_{ \la \in E_n/\sim} \g^{\la} \times
\prod_{\la \in F_n} \g^{\la},
$$
c'est-à-dire que $\g'_n$ contient une algèbre de Lie qui,
d'après la proposition \ref{detgl}, a même dimension
que l'ensemble d'arrivée de $\phi_n$, donc $\phi_n$ est surjectif.
\end{proof}

\subsection{Conclusion de la démonstration}

Remarquons d'abord que, d'après le lemme \ref{minork} et si $n \geq 7$,
un idéal simple $\g^{\la}$ ne peut être simultanément de type
$A_r,B_r,C_r$ ou $D_r$, puisque ces types sont deux-à-deux
non isomorphes pour $r \geq 4$.

Supposant $\HRm$ il s'agit maintenant, d'après le lemme \ref{lemconclu}, de montrer
que la propriété $\mathrm{(IS)_n}$ est vérifiée, c'est-à-dire
que, si $\la, \mu \vdash n$ ne sont pas toutes deux des équerres
et $\g^{\la} = \g^{\mu}$, alors $\mu \in \{ \la, \la' \}$. Soient
donc deux tels $\la, \mu$.

Si $\la$ ou $\mu$ est une équerre, on peut donc supposer que seul
$\la \not\in E_n \cup F_n$. Ainsi $\g^{\la} = \g^{\alpha}$
est de type $A_{n-2}$, donc ne peut
être isomorphe à une algèbre de Lie orthogonale ou
symplectique. On en déduit $\mu \in E_n$, et ainsi
$\g^{\mu}$ est de type $A_{\dim(\mu) -1 }$. Alors
$\g^{\la} = \g^{\mu}$ implique $n-2 = \dim \mu -1$ soit $\dim \mu =
n-1$, ce qui est impossible d'après le lemme \ref{ppdim}
puisque $\mu \in E_n$.

On suppose désormais $\la, \mu \in E_n \cup F_n$. D'après la remarque
précédente on a alors, soit $\la, \mu \in E_n$, soit $\la,\mu \in F_n$.
Les représentations $\la, \la', \mu, \mu'$ de $\g'_n$ se factorisent
dans tous les cas par un même idéal simple $\a = \g^{\la} = \g^{\mu}$
de $\g'_n$.

Si $\la, \mu \in E_n$, alors l'idéal simple $\a$ est de type
$A_{N-1}$, avec $N = \dim(\la) = \dim(\mu)$. On a
$N \geq 5$ d'après le lemme \ref{minork}. Si
$\mu \not\in \{ \la , \la' \}$, ce qui est équivalent
à $\mu' \not\in \{ \la , \la' \}$, alors $\la,\la', \mu$ seraient
trois représentations irréductibles de $\a$ de dimension $N$, deux-à-deux
non isomorphes d'après la proposition \ref{equivisom}. C'est exclu car une
algèbre de Lie simple de type $A_{N-1}$ n'admet à isomorphisme
près que deux représentations irréductibles de dimension $N$, pour tout $N \geq 3$.

Si $\la,\mu \in F_n$, alors $\a$ est de type $B_r,C_r$ ou $D_r$ pour
$r = N/2$, avec $N = \dim(\la) = \dim(\mu)$, d'après la proposition
\ref{detgl}. De plus $r \geq 8$ d'après le lemme \ref{minork}. Or,
pour $r \geq 5$ les algèbres de Lie simples de type $B_r,C_r$
ou $D_r$ n'admettent à isomorphisme près qu'une représentation irréductible
de dimension $2r$. Ainsi $\la$ et $\mu$ sont isomorphes en tant
que représentations de $\g'_n$, donc de $\SN$ d'après
la proposition \ref{equivisom}, soit $\la = \mu$. Ceci conclut
la démonstration du théorème A.

\section{Enveloppes algébriques}

\subsection{Rappels et compléments sur les groupes algébriques}

Soient $K$ un corps de caractéristique 0 et $G$ un
sous-groupe de $GL_N(K)$. On note $\bar{G}$ et
on appelle enveloppe algébrique de $G$ le plus petit sous-groupe
algébrique de $GL_N(K)$ qui contient $G$. Si $H$ est un sous-groupe
(distingué) de $G$, alors $\bar{H}$ est un sous-groupe (distingué) de
$\bar{G}$.

Le groupe $\bar{G}$ est l'adhérence de $G$ dans $GL_N(K)$ pour
la topologie de Zariski. En particulier, il existe des
fonctions polynomiales $a_1,\dots,a_r$ sur $K^{N^2}$ qui définissent
$\bar{G}$, c'est-à-dire telles que
$$ \bar{G} = \{ m \in GL_N(K) \ \mid \ a_1(m) = \dots = a_r(m) = 0 \}
$$
Si $K$ est le corps des fractions d'un anneau $A$, on peut de plus
supposer que $a_1,\dots,a_r$ sont à coefficients dans $A$.

Le groupe algébrique $G$ définit un foncteur $\mathcal{A} \to
\bar{G}(\mathcal{A})$ des $K$-algèbres commutatives à unité vers les
groupes. On appelle $\mathcal{A}$-point de $\bar{G}$ un élément
de $\bar{G}(\mathcal{A})$. On note $\Lie \bar{G}$ l'algèbre de Lie
de $\bar{G}$. Il s'agit d'une $K$-sous-algèbre de Lie de $\gl_N(K)$,
qui peut être définie en caractéristique 0 par la théorie de
l'exponentiation formelle de Chevalley telle qu'établie dans \cite{CHEV} \S 12
$$
\Lie \bar{G} = \{ x \in \gl_N(K) \ \mid \ \exp(ux) \in \bar{G}(K[[u]]) \}
$$

Soit $\k$ un corps, prenons $A = \k[[h]]$ l'anneau des séries
formelles correspondant et $K = \k((h))$ son
corps de fractions. En plus des raisonnement précédents,
on a dans ce cadre le résultat suivant :
\begin{lemme} \label{lemadh} Soit $G$ un sous-groupe de $GL_N(K)$, et
$X = \exp(hx) \in G$ avec $x \in \gl_N(A)$. Alors
$x \in \Lie \bar{G}$.
\end{lemme}
\begin{proof}
Il s'agit de montrer que $\exp(uhx) \in \bar{G}(K[[u]])$, puisque
$hx \in \Lie \bar{G}$ si et seulement si $x \in \Lie \bar{G}$.
Notons $a_1,\dots,a_r$ une famille de fonctions polynomiales définissant
$\bar{G}$, que l'on choisit à coefficients dans $A$. Alors on montre
facilement que $Q_i(u) = a_i(\exp(uhx)) \in M_N(L)$
avec $L = (\k[u])[[h]]$, c'est-à-dire $Q_i(u) = \sum_{j=0}^{\infty}
b_{i,j}(u)h$ avec $b_{i,j} \in \k[u]$. En particulier, $\exp(uhx) \in
\bar{G}(K[[u]])$ si et seulement si $b_{i,j} = 0$ pour tout
$1 \leq i \leq r$ et $j \in \N$. Or, si $X \in G$, on a $X^n
\in G \subset \bar{G}$ pour tout $n \in \Z$, c'est-à-dire
$b_{i,j}(n) = 0$ pour tout $n \in \Z$, ce qui implique $b_{i,j} = 0$
par Zariski-densité de $\Z$ dans tout corps $\k$ de caractéristique 0.
\end{proof}

Autrement dit, on a
$$
\Lie \bar{G} \supset \{ x \in \gl_N(A) \ \mid \ \exp(hx) \in G \}
$$

Nous rappelons également sous forme de lemme le résultat
suivant de la théorie élémentaire des groupes algébriques

\begin{lemme} \label{leminclu} Soient $H$ et $G$ deux sous-groupes algébriques de $GL_N(K)$.
Si $G$ est irréductible, $H \subset G$ et $\Lie H = \Lie G$, alors
$H = G$.
\end{lemme}
\begin{proof}
Soit $H_0$ la composante connexe de $H$. C'est un sous-groupe algébrique
ir\-ré\-duc\-ti\-ble d'indice fini de $H$. On a donc
$\Lie H_0 = \Lie H$ et $\Lie H_0 = \Lie G$ par hypothèse. Comme
$H_0 \subset G$ et que $H_0$ et $G$ sont irréductibles, on en déduit
$H_0 = G$ (cf. \cite{CHEV}, thm. 8 cor. 1), donc $H = G$.
\end{proof}

\subsection{Préliminaires sur les groupes orthogonaux et symplectiques}

\subsubsection{Notations}
Soit $\k$ un corps de caractéristique nulle, et $n = 2p \geq 2$ un
entier positif pair. On munit $\k^n$ de son produit scalaire
standard $<\ , \ >$. On introduit les matrices de 
$GL_n(\k)$ suivantes, définies par blocs de taille $p \times p$ :
$$
D_o = \left( \begin{array}{cc} 1 & 0 \\ 0 & -1 \end{array} \right) \ \ \ 
D_s = \left( \begin{array}{cc} 0 & 1 \\ -1 & 0 \end{array} \right) \ \ \ 
J = \left( \begin{array}{cc} 0 & 1 \\ 1 & 0 \end{array} \right)
$$
et l'on note également $H \in GL_n(\k)$ la matrice
diagonale dont le premier coeffcient vaut $-1$ et les suivants $1$.
Outre $SL_n(\k)$, un autre sous-groupe de $GL_n(\k)$ important
pour nous est
$$
\widetilde{SL}_n(\k)
= \{ A \in GL_n(\k) \ | \ \det A = \pm 1 \}
$$
dont $SL_n(\k)$ est un sous-groupe d'indice 2.

\subsubsection{Cas orthogonal}
Soit $(\ | \ )$ la forme bilinéaire symétrique définie sur $\k^n$ par
$(x|y) = < D_0 x | y>$. Pour tout $A \in GL_n(\k)$, on note $A^*$
son adjoint par rapport à cette forme. Cela nous permet
d'introduire les groupes
$$
O_{p,p}(\k) = \{ A \in GL_n(\k) | A^{-1} = A^* \}\ \ \ 
\widetilde{O}_{p,p}(\k) = \{ A \in GL_n(\k) | A^{-1} = \pm A^* \}
$$
On a un morphisme $\varphi : \widetilde{O}_{p,p}(\k) \to \{ \pm 1 \}$
défini par $\varphi(A) = u$ si $A^{-1} = u A^*$, qui
est surjectif car $J \in \widetilde{O}_{p,p}(\k)$
et $\varphi(J) = -1$. On en déduit une suite exacte scindée
$$
1 \to O_{p,p}(\k) \to \widetilde{O}_{p,p}(\k) \to \{ \pm 1\} \to 1
$$

Le groupe spécial orthogonal est défini par $SO_{p,p}(\k)
= \{ A \in O_{p,p}(\k) \ | \ \det (A) = 1 \}$. L'application
$\varphi \times \det$ induit alors une suite
$$
1 \to SO_{p,p}(\k) \to \widetilde{O}_{p,p}(\k) \to \{ \pm 1 \} ^2 \to 1
$$
qui est exacte et scindée pour tout $p\geq 1$. En effet, on vérifie
facilement
$$
\left\lbrace \begin{array}{lcl} \det(J) & = & (-1)^p  \\
\varphi(J) & = & -1 \end{array} \right.
\left\lbrace \begin{array}{lcl} \det(H) & =& -1  \\
\varphi(H) & = & 1 \end{array} \right.
\left\lbrace \begin{array}{lcl} \det(JH) & = &(-1)^{p+1}  \\
\varphi(JH) & = & -1 \end{array} \right.
$$

A chacun des trois sous-groupes propres de $\{ \pm 1 \}^2$ sont ainsi
associés trois sous-groupes (distingués)
d'indice deux de $\widetilde{O}_{p,p}(\k)$
qui contiennent $SO_{p,p}(\k)$. Ce sont
$$
\begin{array}{lcl}
\widetilde{SO}^1_{p,p}(\k) & = & SL_n(\k) \cap \widetilde{O}_{p,p}(\k) \\
\widetilde{SO}^2_{p,p}(\k) & = & O_{p,p}(\k) \\
\widetilde{SO}^3_{p,p}(\k) & = & \{ A \in
\widetilde{SL}_n(\k) | A^{-1} = \det(A) A^* \} \\
\end{array}
$$

Par convention, on notera ici 
$$
\begin{array}{lcll}
\widetilde{SO}_{p,p}(\k)
& = & \widetilde{SO}^1_{p,p}(\k) & \mbox{ si $p$ est pair,} \\
\widetilde{SO}_{p,p}(\k) & = & \widetilde{SO}^3_{p,p}(\k) & \mbox{ si
$p$ est impair.}
\end{array}
$$

\subsubsection{Cas symplectique}

On considère la forme bilinéaire antisymétrique non dégénérée définie
sur $\k^n$ par $(x|y) = <D_s x ,y>$. On note encore $A^*$
l'adoint de $A$ par rapport à $(\ | \ )$ et on définit
les groupes
$$
\begin{array}{lcl}
Sp_n(\k) & = & \{ A \in GL_n(\k) \ | \ A^{-1} = A^* \} \\
\widetilde{Sp}_n(\k) & = & \{ A \in GL_n(\k) \ | \ A^{-1} = \pm A^* \} 
\end{array}
$$
On a $Sp_n(\k) \subset SL_n(\k)$, 
$\widetilde{Sp}_n(\k) \subset \widetilde{SL}_n(\k)$
et $Sp_n(\k) \neq \widetilde{Sp}_{n}(\k)$ parce que $J \in \widetilde{Sp}_n(\k)
\setminus Sp_n(\k)$. On a ainsi une suite exacte scindée
$$
1 \to Sp_n(\k) \to \widetilde{Sp}_n(\k) \to <J>
$$

\section{Groupes de tresses et algèbres de Hecke}

Dans ce qui suit, on fixe un corps $\k$ de caractéristique
nulle. L'objet de cette section est de rappeler
des propriétés des groupes de tresses et leur relation
avec l'algèbre de Hecke de type $A$. Notamment, on

\subsection{Groupes de tresses}

Pour tout $n \geq 2$, on note $\B_n$ le groupe de tresses à $n$
brins, défini
par générateurs $\sigma_1, \dots, \sigma_{n-1}$ et
relations
$$
\sigma_i \sigma_j = \sigma_j \sigma_i \mbox{ pour }
|i-j| \geq 2 , \mbox{ et } \sigma_i \sigma_{i+1} \sigma_i = \sigma_{i+1}
\sigma_i \sigma_{i+1} \mbox{ pour }1 \leq i \leq n-2.
$$
Une conséquence immédiate de ces relations est que les éléments
$\sigma_i$ pour $1 \leq i \leq n-1$, appelés générateurs d'Artin, sont
conjugués entre eux.
Pour les autres résultats classiques concernant ce groupe, que nous ne
faisons que rappeler, on pourra
se reporter à \cite{refbn}.
Le morphisme
surjectif $\B_n \to \SN$ défini par $\sigma_i \mapsto s_i = (i \ i+1)$
admet pour noyau le groupe de tresses pures, noté $\P_n$,
qui est engendré par les éléments $\xi_{ij}$, pour $1 \leq i < j \leq n$,
définis par
$$
\xi_{ij} = \sigma_{j-1} \dots \sigma_{i+1} \sigma_i^2 \sigma_{i+1}^{-1}
\dots \sigma_{j-1}^{-1}.
$$

A ce groupe est associée une $\k$-algèbre de Lie, dite des tresses
infinitésimales pures, définie par générateurs
$t_{ij}$, $1 \leq i,j \leq n$,
et relations
$$
\begin{array}{l}
t_{ij}  = t_{ji}, t_{ii} = 0,\\
 {[}t_{ij} ,t_{kl}] = 0
\mbox{ si }\#\{i,j,k,l\} = 4, \\
{[}t_{ij},t_{ik} + t_{kj}] = 0 \mbox{  pour tous }i,j,k
\end{array}
$$
Elle est graduée, ainsi que son algèbre enveloppante universelle
$\U \mathcal{T}_n$ par $\deg(t_{ij})=1$.
On note $\widehat{\mathcal{T}_n}$ et $\widehat{\U \mathcal{T}_n}$
les complétions de ces algèbres par
rapport à leurs graduations respectives, $\exp \mathcal{T}_n
\subset \widehat{\U \mathcal{T}_n}$ le groupe pro-algébrique
associé à $\mathcal{T}_n$.

Cette algèbre de Lie, son algèbre enveloppante et
leurs complétions respectives sont de plus
naturellement munies d'une action par automorphismes naturelle de $\SN$, par
$s.t_{i,j} = t_{s(i),s(j)}$ pour tout $s \in \SN$.
Cela permet d'introduire une $\k$-algèbre de Hopf
complétée $\k \SN \ltimes \U \widehat{\mathcal{T}_n}$, dans
laquelle on peut construire le groupe $\SN \ltimes
\exp \widehat{\mathcal{T}_n}$. Dans \cite{DRIN},
V. Drinfel'd a défini des morphismes
injectifs
$\widetilde{\Phi} : \B_n \to \SN \ltimes \exp \widehat{\mathcal{T}_n}$
tels que
$
\PPhi(\sigma_i) = \Phi_i s_i \exp(t_{i,i+1}) \Phi_i^{-1}$ pour tout
$1 \leq i \leq n-1$, avec $\Phi_i \in \exp \TTN$ et $\Phi_1 = 1$. On fixe un tel
morphisme.

On vérifie facilement $\PPhi(\P_n) \subset \exp \TTN$, et
$\PPhi(\xi_{ij}) = \exp \psi_{ij}$ avec $\psi_{ij} \in 2 t_{ij} + 
\widehat{[\mathcal{T}_n,\mathcal{T}_n]}$.
Un sous-groupe distingué remarquable de $\P_n$ est le sous-groupe
$\FN$ engendré par les éléments $\xi_{1,n}, \dots, \xi_{n-1,n}$. C'est un groupe
libre en ces $n-1$ générateurs. Parallèlement, il est facile de
vérifier que la sous-algèbre de Lie $\mathcal{F}_n$ de $\mathcal{T}_n$
engendrée par les éléments $t_{1,n},\dots,t_{n-1,n}$ est
libre en ces générateurs et est un idéal de $\mathcal{T}_n$.
Il s'ensuit que $\PPhi(\FN) \subset \exp \widehat{\mathcal{F}}_n$.

Un dernier sous-groupe qui intervient ici est le sous-groupe $\B_n^2$
des tresses paires. C'est le noyau de la composition du
morphisme naturel $\pi : \B_n \to \SN$ et de la signature
$\SN \to \{ \pm 1 \}$, c'est-à-dire $\pi^{-1}(\ALTN)$
où $\ALTN$ désigne le groupe alterné sur $n$ lettres. On a
évidemment les inclusions
$$
\FN \subset \P_n \subset \B_n^2 \subset \B_n
$$

\subsection{Algèbres de Hecke}

Pour tout $n \geq 3$, on définit habituellement l'algèbre
de Hecke de type $A_{n-1}$, que l'on note $H_n(q)$, sur
un corps $K$ de caractéristique nulle, dans lequel
on a choisit un élément $q$ transcendant sur le corps premier.
C'est une $K$-algèbre définie
par générateurs $\s_1,\dots, \s_{n-1}$ et relations
$$
\begin{array}{ll}
\s_i \s_j = \s_j \s_i & \mbox{ pour }
|i-j| \geq 2 \\
\s_i \s_{i+1} \s_i = \s_{i+1}
\s_i \s_{i+1} & \mbox{ pour }1 \leq i \leq n-2 \\
(\s_i-q)(\s_i+q^{-1}) = 0 & \mbox{ pour }1 \leq i \leq n-1
\end{array}
$$
Il est clair que $H_n(q)$ est un quotient de l'algèbre de groupe
$K \B_n$ de $\B_n$ sur $K$ : plus précisément, c'est le quotient
de $K \B_n$ par l'idéal bilatère engendré par les éléments
$(\sigma_i -q)(\sigma_i + q^{-1})$. Comme les générateurs
d'Artin sont conjugués entre eux, cet idéal est engendré
par un seul élément, par exemple $(\sigma_1 -q)(\sigma_1 + q^{-1})$.

Ainsi, toute représentation
de $H_n(q)$ sur $K$ s'étend en une représentation de $K \B_n$ et réciproquement
toute représentation de $K \B_n$ dans laquelle l'élément
$(\sigma_1 -q)(\sigma_1 + q^{-1})$ agit trivialement se factorise par
$H_n(q)$.

A isomorphisme près, $H_n(q)$ ne dépend pas de l'élément transcendant
$q$ que l'on a choisi, et peut être définie sur $\k(q)$ ; le
choix du corps $K$ est donc de peu d'importance. Pour nos
besoins, il sera commode de prendre pour $K$ le corps
$\k((h))$ des séries de Laurent en une indéterminée sur $\k$,
et $q = e^h$.

Un résultat de Tits (\cite{BOURB} ch. 4 ex. 26 et ex. 27) montre de plus que $H_n(q)$ est de dimension
finie et isomorphe à l'algèbre
de groupe $K \SN$. On peut obtenir un tel isomorphisme à l'aide
du morphisme $\PPhi$ choisi précédemment.

En effet, notant $A = \k[[h]]$ l'anneau des séries formelles en une
indéterminée sur $\k$, on vérifie facilement que l'application $t_{ij} \to h (i \ j)$
s'étend naturellement en un morphisme d'algèbres de Lie surjectif
de $\mathcal{T}_n$ vers $\g_n^h = \g_n \otimes A$, donc de l'algèbre
enveloppante
$\U \mathcal{T}_n$ vers l'algèbre de groupe $A \SN$ de $\SN$
sur l'anneau $A = \k[[h]]$ des séries formelles. Ce morphisme étant continu
pour la topologie $h$-adique respectivement à la graduation
de $\U \mathcal{T}_n$, il s'étend à $\widehat{\U \mathcal{T}}_n$.
Comme ce morphisme est $\SN$-équivariant, $\PPhi$ permet d'en déduire un
morphisme $A \B_n \to A \SN$, dont la
réduction à $h=0$ est surjective.
On en déduit un morphisme surjectif $K \B_n \to K \SN$.
L'image de $\sigma_1$
est alors $X = s_1 \exp(h s_1)$. De $s_1^2 = 1$ on déduit
$$
\exp(h s_1 ) = \frac{q+q^{-1}}{2} + \frac{q-q^{-1}}{2} s_1
$$
et finalement $(X-q)(X+q^{-1}) = 0$. On en déduit un morphisme
surjectif explicite $H_n(q) \to K \SN$, qui est un isomorphisme
par égalité des dimensions, et pour lequel
l'image de $\s_1$ est $s_1 \exp (h s_1)$.

Pour étudier l'image du groupe de tresses dans $H_n(q)$
il suffit donc d'étudier l'image de $\B_n$ par
$\PP : K \B_n \to K \SN$ obtenu par composition du morphisme
précédent avec la projection naturelle $K \B_n \to H_n(q)$.
Puisque $\PPhi(\P_n) \subset \exp \widehat{\mathcal{T}_n}$,
on sait d'ores et déjà que $\PP(\P_n) \subset \exp \g_n^h$.

\subsection{Représentations de l'algèbre de Hecke}

A toute représentation $\rho : \SN \to GL(V)$, où $V$
est un espace vectoriel de dimension finie sur $\k$, on peut associer
par extension des scalaires une représentation encore
notée $\rho$ de $K \SN$ sur $V^h = V \otimes K$.

On en déduit une représentation $R = \rho \circ \PP$ de $\B_n$ sur $V^h$
telle que, si l'on note $x = \rho(s_1)$, alors
$R(\sigma_1) = x e^{hx}$. En particulier,
$
\det(R(\sigma_i)) = \det(x) \exp(h \tr(x))
$ pour tout $1 \leq i \leq n-1$ puisque les
générateurs d'Artin sont conjugués entre eux.

Considérons le cas où $\rho$ est une représentation irréductible,
associée à une partition $\la$ de $n$. Alors $\tr(x) = \gamma(\la)$,
valeur du caractère associé à $\la$ 
sur la classe de conjugaison de $\SN$ formée des
transpositions. Comme $x^2 = 1$, la valeur $\eta(\la) = \det(x)$
appartient à $\{ \pm 1 \}$.

Les valeurs $\dim(\la)$, $\gamma(\la)$ et $\eta(\la)$ sont
reliées entre elles. En effet, si l'on note $a_{\pm}$ la
multiplicité avec laquelle $\pm 1$ intervient dans le spectre
de $x$, on a $\dim(\la) = a_+ + a_-$ et
$\gamma(\la) = a_+ - a_-$, donc
$$
\eta(\la) = (-1)^{a_-} = (-1)^{\frac{\dim(\la) - \gamma(\la)}{2}}
$$

Si $\la = \la'$, les tableaux standards de forme $\la$ et les formes
bilinéaires définies en section 4.2 permettent d'identifier $\la^h$
à $K^N$ pour $N = \dim(\la)$, muni de l'une des formes bilinéaires
étudiée en 7.2. Si cette forme est symérique, on note
$$
OSP(\la^h) = SO_{\frac{N}{2},\frac{N}{2}}(K)   \ \ \ \ \ \ \ \ 
\widetilde{OSP}(\la^h) = \widetilde{SO}_{\frac{N}{2},\frac{N}{2}}(K)
$$
et si elle est antisymétrique,
$$
OSP(\la^h) = Sp_{N}(K)   \ \ \ \ \ \ \ \  
\widetilde{OSP}(\la^h) = \widetilde{Sp}_{N}(K).
$$
On rappelle que dans ces deux cas on a une suite exacte scindée
$$
1 \to OSP(\la^h) \to \widetilde{OSP}(\la^h) \to \{ \pm 1 \} \to 1
$$
Le tableau \ref{tableinc} résume les propriétés que vérifient, dans ces
cas particuliers, l'image de $\B_n$. Leur
démonstration est l'objet des alinéas qui suivent.
\begin{table}
$$
\begin{array}{|c|c|l|c|}
\hline
\gamma(\la) = 0 & \eta(\la) = 1 & R(\B_n) \subset SL(\la^h) & R(\B_n^2) \subset SL(\la^h) \\
\cline{2-3}
 & \eta(\la) = -1 & R(\B_n) \subset \widetilde{SL}(\la^h) \mbox{ et } R(\B_n) \not\subset SL(\la^h) &  \\
\hline
\multicolumn{2}{|c|}{\la = \la'} & R(\B_n) \subset \widetilde{OSP}(\la^h) \mbox{ et } R(\B_n) \not\subset OSP(\la^h) & R(\B_n^2) \subset OSP(\la^h) \\
\hline
\end{array}
$$
\caption{Inclusions remarquables.}
\label{tableinc}
\end{table}

\subsubsection{Cas $\gamma(\la) = 0$.} Si $\gamma(\la) = 0$, alors $\det R(\sigma_i) = \eta(\la)$
pour tout $1 \leq i \leq n-1$, donc $R(\sigma_i) \in \widetilde{SL}(\la^h)$ et donc
$R(\B_n) \subset \widetilde{SL}(\la^h)$. En revanche, si $\eta(\la) = -1$, alors
$R(\sigma_i) \in \widetilde{SL}(\la^h) \setminus SL(\la^h)$.

Dans tous les cas, $\det \rho(a) = 1$ pour tout $a \in \ALTN$, donc $R(\B_n^2) \subset SL(\la^h)$.

\subsubsection{Cas $\la = \la'$.} Si $\la = \la'$, alors $\gamma(\la) = 0$. Cette situation est donc
un cas particulier du paragraphe précédent. De plus $\dim(\la)$ est pair
et, posant
$\dim(\la) = 2p$, on a donc $\eta(\la) = (-1)^p$. Puisque
l'on a $-\rho(s_1)^*
= \rho(s_1) = \rho(s_1)^{-1}$, cela implique par la définition choisie de $\widetilde{OSP}(\la^h)$
que $\rho(s_1) \in \widetilde{OSP}(\la^h) \setminus OSP(\la^h)$. Comme
d'autre part $\rho(\g_n) \subset \osp(\la)$, on en déduit
que $R(\sigma_1) \in
\widetilde{OSP}(\la^h) \setminus OSP(\la^h)$.

Enfin, pour tout $a \in \ALTN$, on a $\det \rho(a) = 1$ et,
chaque permutation paire étant produit d'un nombre pair
de transpositions consécutives, de $- \rho(s_i)^* = \rho(s_i)^{-1}$
on déduit $\rho(a)^* = \rho(a)^{-1}$. Il s'ensuit que
$R(\B_n^2) \subset OSP(\la^h)$.

\subsection{Enveloppe algébrique dans les représentations irréductibles}

L'objet de cette section est, notant $R_{\la}$ la représentation
de $\B_n$ déduite comme pré\-cé\-dem\-ment d'une partition propre $\la$ de
$n$, de
déterminer l'enveloppe algébrique des images de $\B_n$ et de ses
sous-groupes remarquables par $R_{\la}$. Rappelons que
$R_{\la}$ s'identifie à la représentation de $H_n(q)$ associée à $\la$.

\begin{theoB}
Soit $n \geq 3$. Pour toute partition propre $\la \vdash n$, les enveloppes
algébriques de $R_{\la}(\FN)$,
$R_{\la}(\P_n)$ et $R_{\la}(\B_n^2)$
sont égales à un même groupe algébrique connexe $G_{\la}$,
dont l'algèbre
de Lie est l'image de $\g_n$ dans $\gl(\la)$, tensorisée par $K$.
Suivant les valeurs de $\la$,
$G_{\la}$ et l'enveloppe algébrique $\widetilde{G_{\la}}$
de $R_{\la}(\B_n)$ valent
$$
\begin{array}{|c|c|c|c|c|c|}
\hline
 & \gamma(\la) & \eta(\la) & G_{\la} & \widetilde{G_{\la}} & \mbox{ Exemple de tels $\la$} \\
\hline
 \la \neq \la'& \neq 0 &  & \multicolumn{2}{|c|}{GL(\la^h)} & [3,2] \\
\cline{2-6}
   & 0 & 1 &  \multicolumn{2}{|c|}{SL(\la^h)} & [6,3,2,2,2] \\
\cline{3-6}
  & &
-1 & SL(\la^h) & \widetilde{SL}(\la^h) & [9,3,3,3,3,1,1,1] \\
\cline{1-6}
  \la=\la' & & & OSP(\la^h) & \widetilde{OSP}(\la^h) & [2,2] \\
\hline   
\end{array}
$$
\end{theoB}
\begin{proof} On suppose d'abord $\gamma(\la) \neq 0$. Cela
implique $\la \in E_n$, donc $\g_{\la} = \sl(\la)$ d'après
le théorème A. D'autre part l'action de $T_n$ est un scalaire non
nul d'après le lemme \ref{lemactTN}. Comme $\g_{\la} = \sl(\la)$
on en déduit que l'image de $\g_n$ dans $\gl(\la)$ est
égale à $\gl(\la)$. D'après le lemme \ref{lemadh} on en déduit
$\Lie \overline{R(\P_n)} \supset \gl(\la^h)$ donc
$\overline{R(\P_n)} = GL(\la^h)$ d'après le lemme \ref{leminclu}.

On suppose désormais $\gamma(\la) = 0$. Alors l'image de $\g_n$ dans
$\gl(\la)$ est $\g_{\la}$. Deux cas se présentent.

Si $\la \in E_n$, alors $\g_{\la} = \sl(\la)$, donc $\Lie \overline{R(\P_n)}
\supset \sl(\la^h)$. Or $R(\P_n) \subset R(\B_n^2) \subset SL(\la^h)$
donc $\overline{R(\P_n)} = \overline{R(\B_n^2)} = SL(\la^h)$. Si
$\eta(\la) = -1$, alors $R(\B_n) \subset \widetilde{SL}(\la^h)$
et $R(\B_n) \not\subset SL(\la^h)$, donc $SL(\la^h) \subsetneq
\overline{R(\B_n)} \subset \widetilde{SL}(\la^h)$ et
$\overline{R(\B_n)} = \widetilde{SL}(\la^h)$ parce que $SL(\la^h)$
est d'indice 2 dans $\widetilde{SL}(\la^h)$. Si $\eta(\la) = -1$,
on a $R(\B_n) \subset SL(\la^h)$ donc $\overline{R(\B_n)} = SL(\la^h)$.

Si $\la \in F_n$, alors $\g_{\la} = \osp(\la)$, donc
$\Lie R(\P_n) \supset \osp(\la^h)$. Or $R(\B_n^2) \subset OSP(\la^h)$
donc $\overline{R(\P_n)} = OSP(\la^h)$. Puisque $R(\B_n) \subset
\widetilde{OSP}(\la^h)$ mais
$R(\B_n) \not\subset OSP(\la^h)$, on en déduit $OSP(\la^h) \subsetneq
\overline{R(\B_n)} \subset \widetilde{OSP}(\la^h)$ et
$\overline{R(\B_n)} = \widetilde{OSP}(\la^h)$.

Dans tous les cas, les enveloppes algébriques de $R(\P_n)$ et
$R(\B_n^2)$ sont donc égales, et sont des groupes connexes
dont l'algèbre de Lie est l'image de $\g_n$ dans $\gl(\la)$,
tensorisée par $K$. Comme $\g_n$ est engendrée par les
transpositions consécutives et $\overline{R(\FN)}
\subset \overline{R(\P_n)}$
on en déduit $\overline{R(\FN)} = \overline{R(\P_n)} = \overline{R(\B_n^2)}$.
\end{proof}

\section{Enveloppe algébrique dans $H_n(q)$}

Dans cette section, nous décrivons l'enveloppe algébrique
de l'image du groupe de tresses et de ses sous-groupes remarquables
dans l'algèbre de Hecke. Plus précisément nous
décrivons, ce qui est équivalent, l'enveloppe algébrique de l'image
de ces groupes par le morphisme $\PP : K \B_n \to (K \SN)$.

Comme on a un isomorphisme naturel
$$
(K \SN)^{\times} \simeq \prod_{\la \vdash n} GL(\la^h)
$$
on identifiera un élément $x \in (K \SN)^{\times}$
à la collection correspondante $(x_{\la})$ pour $\la \vdash n$,
avec $x_{\la} \in GL(\la^h)$. Nous commençons par
décrire un sous-groupe algébrique de $(K \SN)^{\times}$.

Pour ce faire, on munit comme en 5.1 chaque $\la$ d'un produit scalaire
$\SN$-invariant, que l'on étend en une produit scalaire de $\la^h$. Un tel produit scalaire n'est bien défini
qu'à multiplication par un scalaire près, en revanche l'adjoint
$x_{\la}^{\#}$ d'un élément $x_{\la} \in GL(\la^h)$ est ne dépend
pas du produit scalaire choisi.

On définit $E_n^+ = \{ \la \in E_n | \la < \la'\}$ et $E_n^- = E_n
\setminus E_n^+$. On introduit alors pour tout $\la \in E_n^+$
les morphismes de $\k$-espaces vectoriels $\la' \to \la$
définis en 5.1, que l'on étend en des isomorphismes de $K$-espaces
vectoriels de $(\la')^h$ vers $\la^h$. 

On rappelle d'autre part que, en tant qu'espace vectoriel et $\SN$-module,
$\alpha_r$ s'identifie à la puissance extérieure $\Lambda^r \alpha$.
On en déduit un morphisme $\Lambda^r : GL(\alpha^h) \to GL(\Lambda^r \alpha^h)
= \alpha_r^h$, défini par $\Lambda^r y = y \wedge \dots \wedge y$.
 
\begin{defi} Pour tout $n \geq 3$, on note $G_n(q)$
l'ensemble des $(x_{\la}) \in (K \SN)^{\times}$
qui vérifient
\begin{enumerate}
\item $P_{\la} x_{\la'}^{-1} =  x_{\la}^{\#} P_{\la}$ pour tout $\la \in E_n^+$.
\item $x_{\la} \in OSP(\la^h)$ pour tout $\la \in F_n$.
\item $\det(x_{\la}) = x_{[n]}^{\gamma(\la)}$ pour tout $\la \vdash n$.
\item $x_{[n]}^{r-1} x_{\alpha_r} = \Lambda^r x_{\alpha}$ pour tout $r \in [1,n]$.
\end{enumerate}
\end{defi}

Pour étudier en particulier les représentations correspondant aux équerres,
nous aurons besoin d'un énoncé plus précis que le lemme \ref{isoequerres}.
On introduit, pour tout $a \in \k$, la représentation $\psi_a$ de dimension 1
de l'algèbre de Hopf $\k \SN \ltimes \U \mathcal{T}_n$
définie par
$$
\left\lbrace \begin{array}{lcl}
\psi_a(t_{ij}) = a & \mbox{ pour } & 1 \leq i < j \leq n \\
\psi_a(s) = 1 & \mbox{ pour } & s \in \SN
\end{array} \right.
$$
On remarque que $\psi_a = \psi_1^{\otimes a}$ si $a \in \N$, et
$\psi_1 = \rho_{[n]}$. La démonstration du lemme \ref{isoequerres}
montre en fait le résultat suivant
\begin{lemme} \label{isoequerres2} En tant que représentation de
$\k \SN \ltimes \U \widehat{\mathcal{T}_n}$, on a
$\Lambda^r \alpha = \alpha_r \otimes \psi_{r-1}$.
\end{lemme}
Notant $\Psi_a = \psi_a \circ \PP$, on en déduit immédiatement
que, pour tout $g \in \B_n$,
$$
R_{\Lambda^r \alpha}(g) = R_{\alpha_r}(g) \Psi_{r-1}(g) = R_{\alpha_r}(g)
\Psi_1(g)^{r-1} = R_{\alpha_r}(g) R_{[n]}(g)^{r-1}
$$

\begin{lemme} \label{inclusion} $\PP(\FN) \subset \PP(\P_n) \subset \PP(\B_n^2) \subset G_n(q)$.
\end{lemme}
\begin{proof} Il suffit de montrer $\PP(\B_n^2) \subset G_n(q)$.
Pour tout $\la \vdash n$, notons $\rho_{\la}$ et $R_{\la}$
les représentations de $\g_n$ et $\B_n$ associées.
On a évidemment $\PP(x) = (R_{\la}(x))$ pour tout $x \in \B_n$.

La première équation découle de ce que, notant
$\PP(\Phi_i) = \exp(\varphi_i)$ pour un certain $\varphi_i \in \g_n^h$,
et $\Phi_{i,\mu} = \exp(\rho_{\mu}(\varphi_i))$ on a,
pour tous $i \in [1,n-1]$ et $\mu \vdash n$,
$$
R_{\mu}(\sigma_i) = \Phi_{i,\mu} \rho_{\mu}(s_i) \exp( h \rho_{\mu}(s_i))
\Phi_{i,\mu}^{-1}$$
Or, d'après la section 5.1, pour tout $u \in \g_n$, donc pour tout $u \in \g_n^h$,
on a $P_{\la} \rho_{\la'}(u) P_{\la}^{-1} = -\rho_{\la}(u)^{\#}$.
D'autre part, pour toute transposition $s$ de $\SN$ on a
$P_{\la} \rho_{\la'}(s) P_{\la}^{-1} = -\rho_{\la}(s)$ donc,
pour tout $g \in \ALTN$, on a $P_{\la} \rho_{\la'}(g) P_{\la}^{-1} = 
\rho_{\la}(g)$.

En particulier, pour $\la \in E_n^+$,
$$
((\Phi_{i,\la})^{-1})^{\#} = \exp(-\rho_{\la}(\varphi_i))
= \exp(P_{\la} \rho_{\la'}(\varphi_i) P_{\la}^{-1} ) = P_{\la}
\Phi_{i,\la'} P_{\la}^{-1}
$$
donc $(R_{\la}(\sigma_i)^{-1})^{\#} = -P_{\la} R_{\la'}(\sigma_i)
P_{\la}^{-1}$. On en déduit que, pour tout $g \in \B_n^2$,
$$
R_{\la'}(g)^{-1} = P_{\la}^{-1} R_{\la}(g)^{\#} P_{\la}
$$
c'est-à-dire que $\PP(\B_n^2)$ vérifie la première équation de définition
de $G_n(q)$.

Le
fait que $R_{\la}(\B_n^2) \subset OSP(\la^h)$ pour tout $\la$
tel que $\la = \la'$, donc en particulier pour tout $\la \in F_n$,
a déjà été démontré.

Soit maintenant $g \in B_n^2$. On a $\PP(g) = \pi(g) \exp(t)$
avec $\pi(g) \in \AN$ et $t \in \g_n^h$. Le centre de $\g_n^h$
est de dimension 1, engendré par
$$
\tilde{T} = \frac{2}{n-1} T_n = \frac{2}{n-1} \sum_{1 \leq i < j \leq n}
(i \ j)
$$
D'après le lemme \ref{lemactTN}, on a $\tr \rho_{\la} ( \tilde{T} )
= \gamma(\la)$ et en particulier $\tr \rho_{[n]}(\tilde{T}) = 1$. Comme
$\g_n^h$ est réductive, il existe $m \in K$ et $t' \in [\g_n^h,\g_n^h]$
tels que $t = m \tilde{T} + t'$. En
particulier, $R_{[n]}(g) = \exp(m)$. De plus
$\det \rho_{\la}(\ALTN) = \{ 1 \}$, donc
$$
\det R_{\la}(g) = \exp( m \tr \rho_{\la}(\tilde{T})  ) = \exp(m \gamma(\la))
= R_{[n]}(g)^{\gamma(\la)}
$$
pour tout $g \in \B_n^2$. Enfin, la dernière équation est une
traduction immédiate du fait que $R_{\Lambda^r \alpha} (g) =
R_{\alpha_r}(g) R_{[n]}(g)^{r-1}$ pour tout $g \in \B_n^2$.
\end{proof}

\begin{lemme} \label{connexite} $G_n(q)$ est un sous-groupe algébrique connexe
de $(K \SN)^{\times}$, d'algèbre de Lie $\mathcal{L}_n \otimes K \simeq \g_n^h$.
\end{lemme}
\begin{proof} Comme les équations qui le définissent sont
polynomiales, $G_n(q)$ est fermé dans $(K \SN)^{\times}$, qui
est d'évidence un groupe algébrique connexe. Il est
immédiat qu'il s'agit d'un sous-groupe de $(K \SN)^{\times}$,
donc d'un sous-groupe algébrique de $(K \SN)^{\times}$.
On introduit le groupe algébrique
$$
\tilde{G}_n(q) = \left( \prod_{\la \in E_n^+} SL(\la^h) \right)
\times \left( \prod_{\la \in F_n} OSP(\la^h) \right) \times SL(\alpha^h)
\times K^{\times}
$$
Ce groupe est connexe, d'algèbre de Lie $\mathcal{L}_n \otimes
K$. On a un morphisme naturel $\tilde{G}_n(q) \to G_n(q)$ de groupes
algébriques, donné par
$$
\left( (y_{\la})_{\la \in E_n^+} ,(z_{\la})_{\la \in F_n}, w,u \right) \mapsto
(x_{\la})
$$
avec
$$
x_{\la} = \left\lbrace 
\begin{array}{lcl}
y_{\la} & \mbox{ si } & \la \in E_n^+ \\
   P_{\la}^{-1}((y_{\la'})^{\#})^{-1}P_{\la} & \mbox{ si } & \la \in E_n^- \\
  z_{\la} & \mbox{ si } & \la \in F_n 
\end{array} \right.
$$
et $x_{\alpha_r} = (\Lambda^r w) u^{1-r}$. Il est clair que
cette application est un isomorphisme de groupes algébriques,
dont la réciproque est simplement donnée par 
$$
(x_{\la}) \mapsto \left( (x_{\la})_{\la \in E_n^+},
(x_{\la})_{\la \in F_n}, x_{\alpha}, x_{[n]} \right)
$$
d'où l'on déduit que $G_n(q)$ est connexe d'algèbre de Lie $\mathcal{L}_n
\otimes K$.
\end{proof}

\begin{theoC}
Pour tout $n \geq 3$, les images des sous-groupes $\FN$, $\P_n$ et $\B_n^2$
du groupe de tresses $\B_n$ dans $H_n(q)$ ont même enveloppe
algébrique $G_n(q)$. Ce groupe algébrique est connexe, d'algèbre
de Lie $\g_n^h$, et est d'indice 2 dans l'enveloppe algébrique de
l'image de $\B_n$.
\end{theoC}
\begin{proof}
D'après les deux lemmes précédents,
l'adhérence Zariski $\overline{\PP(\FN)}$ de $\PP(\FN)$ est un
sous-groupe algébrique du groupe algébrique connexe $G_n(q)$.
On a d'autre part
$$
\g_n^h \subset \Lie \overline{\PP(\FN)} \subset \Lie G_n(q) = \mathcal{L}_n
\otimes K.
$$
En effet, $\overline{\PP(\FN)} \subset G_n(q)$ d'après le lemme
\ref{inclusion}, et $\Lie \overline{\PP(\FN)}$ contient l'image
de $\mathcal{F}_n$ par le morphisme naturel $\mathcal{T}_n \to
\g_n$, tensorisée par $K$. Comme, d'après le lemme \ref{gentransp},
cette image est $\g_n$, on a bien $\g_n^h \subset \Lie \overline{\PP(\FN)}$.
Enfin, $\Lie G_n(q) = \mathcal{L}_n
\otimes K$ découle du lemme \ref{connexite}.

On a donc $\Lie \overline{\PP(\FN)} =
\Lie G_n(q)$ puisque $\g_n^h \simeq \mathcal{L}_n \otimes K$ d'après
le théorème A. On déduit alors du lemme \ref{leminclu} et de la connexité de
$G_n(q)$ que $\overline{\PP(\FN)} = G_n(q)$, donc
$\overline{\PP(\FN)} = \overline{\PP(\P_n)} = \overline{\PP(\B_n^2)}
= G_n(q)$ d'après le lemme \ref{inclusion}.

On a $\B_n = \B_n^2 \sqcup \sigma_1 \B_n^2$, et
$\PP(\sigma_1) \not\in \overline{\PP(\B_n^2)}$ car pour tout $n \geq 3$
il existe $\la_0 \vdash n$ tel que $\la_0 = \la_0'$. En effet,
si $n = 2p+1$ est impair on peut choisir $\la_0 = [p+1,p]$ et
si $n=2p$ est pair avec $p \geq 2$ on peut choisir $\la_0 = [p,2,p-2]$.
On en déduit $R_{\la_0}(\sigma_1) \not\in OSP(\la^h)$ d'après le tableau
\ref{tableinc},
alors que $\overline{R_{\la_0}(\B_n^2)}$, projection de
$\overline{\PP(\B_n^2)}$ sur $GL(\la_0^h)$, est incluse dans $OSP(\la^h)$.

Introduisant le groupe algébrique $\tilde{G}_n(q) = G_n(q) \sqcup
\PP(\sigma_1) G_n(q)$ on en déduit que $\overline{\PP(\B_n)}=
\tilde{G}_n(q) \neq G_n(q)$, et que $\overline{\PP(\B_n)}$
contient bien $G_n(q)$ comme sous-groupe connexe d'indice 2.
\end{proof}

\end{document}